\def\coralreport{1}

\documentclass[10pt]{article}

\usepackage{fec,hyperref,ifthen}

\ifthenelse{\coralreport = 1}{
  \usepackage{isetechreport}
  \def\reportyear{17T}
  \def\reportno{17}
  \def\revisionno{1}
  \def\originaldate{January 27, 2018}
  \def\revisiondate{\today}
  \coraltrue
  \cvcrfalse
  \RequirePackage{amsthm}
  \hypersetup{colorlinks,urlcolor=blue}
  \newtheorem{theorem}{Theorem}[section]

  \numberwithin{equation}{section}
}{
\OneAndAHalfSpacedXI

\usepackage{natbib}
 \bibpunct[, ]{(}{)}{,}{a}{}{,}%
\TheoremsNumberedThrough     
\ECRepeatTheorems

\EquationsNumberedThrough    

\MANUSCRIPTNO{}
}

\usepackage{hyperref}
\usepackage{xcolor}
\hypersetup{colorlinks}
\usepackage{fec}
\usepackage{subcaption}
\usepackage{tikz}
\usetikzlibrary{calc}

\newcommand{\trish}{\texttt{TRish}}

\begin{document}

\ifthenelse{\coralreport = 1}{

  \title{A Stochastic Trust Region Algorithm Based on Careful Step Normalization}

  \author{Frank E.~Curtis\thanks{E-mail: \texttt{\href{mailto:frank.e.curtis@lehigh.edu}{frank.e.curtis@lehigh.edu}}}}
  \author{Katya Scheinberg\thanks{E-mail: \texttt{\href{mailto:katyas@lehigh.edu}{katyas@lehigh.edu}}}}
  \author{Rui Shi\thanks{E-mail: \texttt{\href{mailto:rus415@lehigh.edu}{rus415@lehigh.edu}}}}
  \affil{Department of Industrial and Systems Engineering, Lehigh University}
  \titlepage

}
{


\RUNAUTHOR{Curtis, Scheinberg, and Shi}

\RUNTITLE{A Stochastic Trust Region Algorithm Based on Careful Step Normalization}

\TITLE{A Stochastic Trust Region Algorithm\\ Based on Careful Step Normalization}

\ARTICLEAUTHORS{
\AUTHOR{Frank E.~Curtis, Katya Scheinberg, and Rui Shi}
\AFF{Department of Industrial and Systems Engineering, Lehigh University, Bethlehem, PA 18015, USA \\ \EMAIL{frank.e.curtis@lehigh.edu}, \EMAIL{katyas@lehigh.edu}, \EMAIL{rus415@lehigh.edu}}
}

}

\maketitle

\ifthenelse{\coralreport = 1}{
\begin{abstract}
An algorithm is proposed for solving stochastic and finite sum minimization problems.  Based on a trust region methodology, the algorithm employs normalized steps, at least as long as the norms of the stochastic gradient estimates are within a specified interval.  The complete algorithm---which dynamically chooses whether or not to employ normalized steps---is proved to have convergence guarantees that are similar to those possessed by a traditional stochastic gradient approach under various sets of conditions related to the accuracy of the stochastic gradient estimates and choice of stepsize sequence.  The results of numerical experiments are presented when the method is employed to minimize convex and nonconvex machine learning test problems.  These results illustrate that the method can outperform a traditional stochastic gradient approach.
\end{abstract}
}{
\ABSTRACT{An algorithm is proposed for solving stochastic and finite sum minimization problems.  Based on a trust region methodology, the algorithm employs normalized steps, at least as long as the norms of the stochastic gradient estimates are within a specified interval.  The complete algorithm---which dynamically chooses whether or not to employ normalized steps---is proved to have convergence guarantees that are similar to those possessed by a traditional stochastic gradient approach under various sets of conditions related to the accuracy of the stochastic gradient estimates and choice of stepsize sequence.  The results of numerical experiments are presented when the method is employed to minimize convex and nonconvex machine learning test problems.  These results illustrate that the method can outperform a traditional stochastic gradient approach.}

\KEYWORDS{stochastic optimization, finite sum minimization, stochastic gradient method, trust region method, machine learning, logistic regression, deep neural networks}


}

\maketitle

\section{Introduction}

The stochastic gradient (SG) method is the signature strategy for solving stochastic and finite-sum minimization problems.  In this iterative approach, each step to update the solution estimate is obtained by taking a negative multiple of an unbiased gradient estimate.  With careful choices for the stepsize sequence, the SG method possesses convergence guarantees and has been employed to great success for solving various types of problems, such as those arising in machine learning.  For fundamental work on SG, see \cite{RobbMonr51} and \cite{RobbSieg71}.

One disadvantage of the SG method is that stochastic gradients, like the gradients that they approximate, possess \emph{no natural scaling}.  By this, we mean that in order to guarantee convergence, the algorithm needs to choose stepsizes in a problem-dependent manner; e.g., common theoretical guarantees require that the stepsize is proportional to $1/L$, where $L$ is a Lipschitz constant for the gradient of the objective function.  This is in contrast to Newton's method for minimization, for which one can obtain (local) convergence guarantees with a stepsize of 1.  Admittedly, Newton's method is not generally guaranteed to converge from remote starting points with unit stepsizes, but these observations do highlight a shortcoming of first-order methods, namely, that for convergence guarantees the stepsizes need always be chosen in a problem-dependent manner.

The purpose of this paper is to propose a new algorithm for stochastic and finite-sum minimization.  Our proposed approach can be viewed as a modification of the SG method.  The approach does not completely overcome the issue of requiring problem-dependent stepsizes, but we contend that our approach does, for practical purposes, reduce somewhat this dependence.  This is achieved by employing, under certain conditions, \emph{normalized} steps.  We motivate our proposed approach by illustrating how it can be derived from a trust region methodology.  This work can be viewed as a first step toward designing new classes of first- and second-order trust region methods for solving stochastic and finite-sum minimization problems.

The use of normalized steps has previously been proposed in the context of (stochastic) gradient methods for solving minimization problems.  For example, in a method that is similar to ours, \cite{hazan2015beyond} propose an approach that employs normalized steps in every iteration.  They show that, if the objective function is \emph{$M$-bounded} and \emph{strictly-locally-quasi-convex}, the stochastic gradients are sufficiently accurate with respect to the true gradients (specifically, when mini-batch sizes are $\Omega(\epsilon^{-2})$), and a sufficiently large number of iterations are run (specifically, $\Omega(\epsilon^{-2})$), then their method will, with high probability, yield a solution estimate that is $\epsilon$-optimal.  By contrast, our approach, by employing a modified update that does not always involve the use of a normalized step, enjoys convergence guarantees under different assumptions.  We argue in this paper that employing normalized steps in all iterations cannot lead to general convergence guarantees, which perhaps explains the additional assumptions required for convergence by \cite{hazan2015beyond}.

It is also worthwhile to mention the broader literature.  For important work on SG-type methods and their corresponding theoretical analyses, see, e.g., \cite{AgarBott15}, \cite{ByrdChinNoceWu12}, \cite{Chun54},  \cite{FrieSchm12}, \cite{GhadLan13}, \cite{Glad65}, \cite{johnson2013accelerating}, and \cite{NemiJudiLanShap09}.  There are also numerous variants of SG methods based on gradient aggregation, iterative averaging, second-order techniques, momentum, acceleration, and beyond; for work on these, see \cite{BottCurtNoce18} and the references therein.  More related to our work are techniques that normalize steplengths based on \emph{accumulated} gradient information; see, e.g., \cite{DuchHazaSing11} and \cite{RossMineLang13}.  In a different direction, one should also contrast our work with stochastic trust region approaches, such as those in \cite{larson2016stochastic} and \cite{chen2015stochastic}.  The approaches proposed in these papers, which are based on the use of randomized models of the objective function constructed during each iteration, are quite distinct from our proposed method.  For example, these approaches follow a traditional trust region strategy of accepting or rejecting each step based on the magnitude of an (approximate) \emph{actual-to-predicted reduction ratio}.  Our method, on the other hand, is closer to the SG method in that it accepts the computed step in every iteration.  Another distinction is that these other approaches rely on the use of so-called \emph{fully linear} models of the objective function to obtain their convergence guarantees.  Our convergence guarantees are obtained under straightforward upper bounds on the second moment of the stochastic gradient estimates, and do not require fully linear models.

The paper is organized as follows.  Our algorithm and motivation for our specific iterate updating scheme are the subject of~\S\ref{sec.algorithm}.  In~\S\ref{sec.convergence}, we prove convergence guarantees for the algorithm under various types of assumptions on the stochastic gradient estimates and stepsize choices.  The results of numerical experiments on test problems---some convex and some nonconvex---are given in~\S\ref{sec.numerical}.  Concluding remarks are given in~\S\ref{sec.conclusion}.  All norms in the paper are Euclidean, i.e., $\|\cdot\| := \|\cdot\|_2$.

\section{Algorithm}\label{sec.algorithm}

Our problem of interest is a stochastic optimization problem in which the goal is to minimize over a vector of decision variables, indicated by $x \in \R{n}$, a function $f : \R{n} \to \R{}$ defined by the expectation of another function $F : \R{n} \times \Xi \to \R{}$ that depends on a random variable $\xi$, i.e.,
\bequation\label{prob.f}
  \min_{x\in\R{n}}\ f(x)\ \ \text{with}\ \ f(x) = \E_\xi[F(x,\xi)],
\eequation
where $\E_\xi[\cdot]$ denotes expectation with respect to the distribution of~$\xi$.  Our algorithm is also applicable for finite-sum minimization where the objective takes the form
\bequation\label{prob.f_sum}
  f(x) = \frac1N \sum_{i=1}^N f_i(x).
\eequation
Such objectives often arise in sample average approximations of~\eqref{prob.f}; e.g., see \cite{ShapDentRusz09}.

\subsection{Algorithm Description}  Our algorithm is stated below as \trish, a trust-region-\emph{ish} algorithm for stochastic optimization.  Each iteration involves taking a step along the negative of a stochastic gradient direction.  In the context of problem~\eqref{prob.f}, this stochastic gradient can be viewed as $g_k = \nabla_x F(x_k,\xi_k)$, where $x_k$ is the current iterate and $\xi_k$ is a realization of the random variable $\xi$.  In the context of problem~\eqref{prob.f_sum}, it can be viewed as $g_k = \nabla_x f_{i_k}(x_k)$ where $i_k$ has been chosen randomly as an index in $\{1,\dots,N\}$.  In addition, in either case, $g_k$ could represent an average of such quantities, i.e., over a set of independently generated realizations $\{\xi_{k,j}\}_{j\in\Scal_k}$ or over independently generated indices $\{i_{k,j}\}_{j\in\Scal_k}$.  This leads to a so-called \emph{mini-batch} approach with $\Scal_k$ representing the mini-batch of samples in the $k$th iteration.  In the algorithm, we simply write $g_k \approx \nabla f(x_k)$ to cover all of these situations, since in any case $g_k$ represents a stochastic gradient estimate for $f$ at $x_k$.

\begin{algorithm}[ht]
  \renewcommand{\thealgorithm}{TRish}
  \caption{(Trust-region-ish algorithm based on careful step normalization)}
  \label{alg.trish}
  \begin{algorithmic}[1]
    \State Choose an initial iterate $x_1$ and positive stepsizes $\{\alpha_k\}$.
    \State Choose positive constants $\{\gamma_{1,k}\}$ and $\{\gamma_{2,k}\}$ such that $\gamma_{1,k} > \gamma_{2,k} > 0$ for all $k \in \N{}$.
    \For{\textbf{all} $k \in \N{} := \{1,2,\dots\}$}
      \State Generate a stochastic gradient $g_k \approx \nabla f(x_k)$.
      \State Set\label{step.cases}
      \bequationn
        x_{k+1} \gets x_k - \left\{ \baligned & \gamma_{1,k} \alpha_k g_k && \text{if $\|g_k\| \in [0,\tfrac{1}{\gamma_{1,k}})$} \\ & \alpha_k g_k/\|g_k\| && \text{if $\|g_k\| \in [\tfrac{1}{\gamma_{1,k}},\tfrac{1}{\gamma_{2,k}}]$} \\ & \gamma_{2,k} \alpha_k g_k && \text{if $\|g_k\| \in (\tfrac{1}{\gamma_{2,k}},\infty)$.} \ealigned \right.
      \eequationn
    \EndFor
  \end{algorithmic}
\end{algorithm}

The scaling of the stochastic gradient employed in \ref{alg.trish} can be motivated in the following manner.  Given a stochastic gradient $g_k$ and a stepsize $\alpha_k$, consider the trust region subproblem
\bequation\label{prob.tr}
  \min_{s\in\R{n}}\ f(x_k) + g_k^Ts\ \ \st\ \ \|s\| \leq \alpha_k.
\eequation
The solution of this subproblem, namely, $s_k = -\alpha_kg_k/\|g_k\|$, represents the step that minimizes the first-order model $f(x_k) + g_k^Ts$ of the objective function $f$ at $x_k$ subject to $s$ having norm less than or equal to $\alpha_k$.  This is the prototypical strategy in a trust region methodology.  When the norm of~$g_k$ falls within the interval $[\frac{1}{\gamma_{1,k}},\frac{1}{\gamma_{2,k}}]$, \ref{alg.trish} takes the step $s_k$.  However, if this were to be done no matter the norm of $g_k$, then the resulting algorithm might fail to make progress in expectation.  This is illustrated in the following example.

\bexample\label{eq.counterexample}
  \textit{
  Suppose that, at a point $x_k \in \R{}$, one has $\nabla f(x_k) = 1$ and obtains
  \bequationn
    g_k = \left\{ \baligned 6 & \text{ with probability $\tfrac13$} \\ -\tfrac32 & \text{ with probability $\tfrac23$}. \ealigned \right.
  \eequationn
  Then, $\E_k[g_k] = 1 = \nabla f(x_k)$, where $\E_k$ denotes expectation given that an algorithm has reached $x_k$ as the $k$th iterate.  However, this means that the normalized stochastic gradient satisfies
  \bequationn
    \frac{g_k}{\|g_k\|} = \left\{ \baligned 1 & \text{ with probability $\tfrac13$} \\ -1 & \text{ with probability $\tfrac23$,} \ealigned \right.
  \eequationn
  from which it follows that $s_k = -\alpha_k g_k/\|g_k\|$ is twice as likely to be a direction of ascent for~$f$ at~$x_k$ than it is to be a direction of descent for $f$ at $x_k$.
  }
\eexample

One can argue from this example that, without potentially restrictive assumptions on the objective function~$f$ and/or the manner in which the stochastic gradient is computed, one cannot expect to be able to prove convergence guarantees for an algorithm that solely computes steps based on solving the trust region subproblem~\eqref{prob.tr}.  In particular, the existence of any point (let alone more than one) at which the expectation is to follow an ascent direction foils the typical convergence theory for an SG approach; see, e.g., \cite{BottCurtNoce18}.

In \ref{alg.trish}, we overcome the issue highlighted in Example~\ref{eq.counterexample} by only choosing the trust region step when the norm of the gradient is within a specified interval; otherwise, we compute a stochastic gradient step with a stepsize that is a multiple of $\alpha_k$.  It is for this reason that we refer to the algorithm as a trust-region-\emph{ish} approach.  Overall, as a function of the norm of the stochastic gradient, the norm of the step taken by the algorithm is illustrated in Figure~\ref{fig.norm_step_vs_norm_g}.  Note that care has been taken to make sure that the norm of the step is a continuous function of the norm of the stochastic gradient estimate.  The plot in Figure~\ref{fig.norm_step_vs_norm_g} illustrates the relationship for moderate values of $(\gamma_{1,k},\gamma_{2,k})$, but notice that for more extreme values (i.e., $\gamma_{1,k} \gg 0$ and $\gamma_{2,k} \approx 0$) the function would essentially be flat (except for stochastic gradients that are very small or large in norm), meaning that the stepsize would typically be scaled so that the step norm is approximately $\alpha_k$ for all $k \in \N{}$.

\bfigure[ht]
  \centering
  \begin{tikzpicture}[scale=2]
  \coordinate (orig) at ( 0.0, 0.0);
  \coordinate [label=right:{$\|g_k\|_2$}] (xmax) at ( 4.0, 0.0);
  \coordinate [label=above:{$\|x_{k+1}-x_k\|_2$}] (ymax) at ( 0.0, 1.5);
  \coordinate (xmin) at (-0.2, 0.0);
  \coordinate (ymin) at ( 0.0,-0.2);
  \draw[thick,->] (xmin) -- (xmax);
  \draw[thick,->] (ymin) -- (ymax);
  \coordinate [label=below:{$\tfrac{1}{\gamma_{1,k}}$}] (g1) at (1.0,0.0);
  \coordinate [label=below:{$\tfrac{1}{\gamma_{2,k}}$}] (g2) at (3.0,0.0);
  \coordinate [label=left:{$\alpha_k$}] (a) at (0.0,1.0);
  \coordinate (b1) at (1.0,1.0);
  \coordinate (b2) at (3.0,1.0);
  \coordinate (b3) at (4.0,1.3333);
  \draw[black!80,thick,->] (orig) -- (b1) -- (b2) -- (b3);
  \draw[black!50,thick,dashed] (a) -- (b1);
  \draw[black!50,thick,dashed] (g1) -- (b1);
  \draw[black!50,thick,dashed] (g2) -- (b2);
  \draw[black!50,thick,dashed] (orig) -- (b2);
  \coordinate [label=above:{$\gamma_{1,k}\alpha_k\ \ $}] (s1) at (0.45,0.55);
  \draw[black!50,thick] (0.45,0.45) -- (s1) -- (0.55,0.55);
  \coordinate [label=above:{$\gamma_{2,k}\alpha_k\ \ $}] (s2) at (3.3,1.2);
  \draw[black!50,thick] (3.3,1.1) -- (s2) -- (3.6,1.2);
\end{tikzpicture}
  \caption{Relationship between $\|g_k\|$ and $\|x_{k+1} - x_k\|$ in Algorithm~\ref{alg.trish}.}
  \label{fig.norm_step_vs_norm_g}
\efigure

Our convergence analysis in the next section requires certain restrictions on the choice of stepsizes---as is typical for (stochastic) gradient methods---and require certains restrictions on $\{\gamma_{1,k}\}$ and $\{\gamma_{2,k}\}$.  For example, the issue in Example~\ref{eq.counterexample} is avoided as long as the pair $(\gamma_{1,k},\gamma_{2,k})$ is chosen such that the step is not normalized with probability 1 at the given $x_k$, which means that---for this particular function, iterate, and variance in the stochastic gradient estimates---one cannot choose this pair such that $|6| \in [\frac{1}{\gamma_{1,k}},\frac{1}{\gamma_{2,k}}]$ and $|\tfrac32| \in [\frac{1}{\gamma_{1,k}},\frac{1}{\gamma_{2,k}}]$ simultaneously.  (In our convergence theory, this is avoided through upper bounds on the ratio $\tfrac{\gamma_{1,k}}{\gamma_{2,k}}$.)  Various situations can illustrate how \ref{alg.trish} avoids the issue in Example~\ref{eq.counterexample}.  For example, consider $\gamma_{1,k} = 1$ and $\gamma_{2,k} = \thalf$, which leads to $\E_k[s_k] = -\thalf\alpha_k(6)(\tfrac13) - \alpha_k(-1)(\tfrac23) = -\tfrac13\alpha_k$, meaning that $s_k$ is a descent direction in expectation.  As another example, consider $\gamma_{1,k} = \tfrac14$ and $\gamma_{2,k} = \tfrac16$, which leads to $\E_k[s_k] = -\alpha_k(1)(\tfrac13) - \tfrac16\alpha_k(-\tfrac32)(\tfrac23) = -\tfrac16\alpha_k$, meaning again that $s_k$ is a descent direction in expectation.  Our theory reveals generic conditions that $\{(\gamma_{1,k},\gamma_{2,k})\}$ must satisfy to attain different convergence properties for \ref{alg.trish}.  We also discuss, in \S\ref{sec.numerical}, strategies for choosing these values in practice.

\section{Convergence Analysis}\label{sec.convergence}

Our goal in this section is to prove convergence guarantees for \ref{alg.trish} that are similar to fundamental guarantees for a straightforward SG method; see, e.g., \cite{BottCurtNoce18}.  As in the notation for Example~\ref{eq.counterexample}, our analysis uses $\E_k[\cdot]$ (resp.~$\P_k[\cdot]$) to denote conditional expectation (resp.~conditional probability) given that the algorithm has reached $x_k$ as the $k$th iterate.

Throughout our analysis, we make the following assumption about the objective function.

\bassumption\label{ass.f}
  The objective $f : \R{n} \to \R{}$ is continuously differentiable and bounded below by $f_* = \inf_{x\in\R{n}} f(x) \in \R{}$.  In addition, at any $x \in \R{n}$, the objective is bounded above by a first-order Taylor series approximation of $f$ at $x$ plus a quadratic term with constant $L \in (0,\infty)$, i.e.,
  \bequation\label{eq.g_Lip}
    f(x) \leq f(\xbar) + \nabla f(\xbar)^T(x - \xbar) + \thalf L \|x - \xbar\|^2\ \ \text{for all}\ \ (x,\xbar) \in \R{n} \times \R{n}.
  \eequation
  It is known that \eqref{eq.g_Lip} holds if the gradient function $\nabla f$ is Lipschitz continuous with constant $L$.  This is often referred to as $L$-smoothness of the function $f$.
\eassumption

We also make the following assumption about the stochastic gradients computed in~\ref{alg.trish}.  This assumption is standard in analyses of SG methods; it is easily seen to be satisfied when the variance of the stochastic gradient estimate is uniformly bounded over $k \in \N{}$.

\bassumption\label{ass.f_limits}
  For all $k \in \N{}$, the stochastic gradient $g_k$ is an unbiased estimator of $\nabla f(x_k)$ in the sense that $\E_k[g_k] = \nabla f(x_k)$.  In addition, there exists a pair $(M_1,M_2) \in (0,\infty) \times (0,\infty)$ (independent of $k$) such that, for all $k \in \N{}$, the squared norm of $g_k$ satisfies
  \bequation\label{eq.second_moment}
    \E_k[\|g_k\|^2] \leq M_1 + M_2 \|\nabla f(x_k) \|^2.
  \eequation
\eassumption

Under these assumptions, we prove the following lemma providing fundamental inequalities satisfied by \ref{alg.trish}.  For ease of reference in this result and throughout the remainder of our analysis, we define the following cases based on those indicated in Line~\ref{step.cases} of \ref{alg.trish}:
\bequationn
  \text{``case 1''}:\ \ \|g_k\| \in [0,\tfrac{1}{\gamma_{1,k}});\ \ \text{``case 2''}:\ \ \|g_k\| \in [\tfrac{1}{\gamma_{1,k}},\tfrac{1}{\gamma_{2,k}}];\ \ \text{``case 3''}:\ \ \|g_k\| \in (\tfrac{1}{\gamma_{2,k}},\infty).
\eequationn
The following lemma reveals an upper bound for the expected decrease in $f$ for all $k \in \N{}$.

\blemma\label{lem.basic}
  Under Assumptions~\ref{ass.f} and \ref{ass.f_limits}, the iterates of \ref{alg.trish} satisfy, for all $k \in \N{}$,
  \bequation\label{eq.cases}
    \baligned
      \E_k[f(x_{k+1})] - f(x_k)
        \leq&\ - \gamma_{1,k} \alpha_k(1 - \thalf \gamma_{1,k} L M_2 \alpha_k) \|\nabla f(x_k)\|^2 \\
            &\ + (\gamma_{1,k} - \gamma_{2,k}) \alpha_k \P_k[E_k] \E_k[\nabla f(x_k)^Tg_k | E_k] + \thalf \gamma_{1,k}^2 L M_1 \alpha_k^2,
    \ealigned
  \eequation
  where $E_k$ is the event that $\nabla f(x_k)^Tg_k \geq 0$ and $\P_k[E_k]$ is the probability of this event.
\elemma
\proof{Proof.}
  For all $k \in \N{}$, let $s_k := x_{k+1} - x_k$ represent the step taken by the algorithm.  By \eqref{eq.g_Lip},
  \bequationn
    f(x_{k+1}) = f(x_k + s_k) \leq f(x_k) + \nabla f(x_k)^T s_k + \thalf L \|s_k\|^2.
  \eequationn
  Thus, with $C_{i,k}$ for $i \in \{1,2,3\}$ respectively representing the events that case 1, case 2, and case 3 occur, and with $\P_k[C_{i,k}]$ for $i \in \{1,2,3\}$ respectively representing the probabilities of these events, one finds from the law of total probability that
  \begin{align}
    \E_k[f(x_{k+1})] - f(x_k)
      &\leq \E_k[\nabla f(x_k)^Ts_k] + \thalf L \E_k[\|s_k\|^2] \nonumber \\
      &=    \sum_{i=1}^3 \P_k[C_{i,k}]\E_k[\nabla f(x_k)^Ts_k | C_{i,k}] + \thalf L \sum_{i=1}^3 \P_k[C_{i,k}]\E_k[\|s_k\|^2 | C_{i,k}]. \label{eq.trick1}
  \end{align}
  In the event $C_{1,k}$, the algorithm yields $s_k = -\gamma_{1,k} \alpha_k g_k$, from which it follows that
  \begin{align}
        &\ \E_k[\nabla f(x_k)^Ts_k | C_{1,k}] \nonumber \\
       =&\ -\gamma_{1,k}\alpha_k\E_k[\nabla f(x_k)^Tg_k | C_{1,k}] \nonumber \\
       =&\ -\gamma_{1,k}\alpha_k\P_k[E_k | C_{1,k}] \E_k[\nabla f(x_k)^Tg_k | C_{1,k} \cap E_k] - \gamma_{1,k} \alpha_k \P_k[\overline{E}_k | C_{1,k}] \E_k[\nabla f(x_k)^Tg_k | C_{1,k} \cap \overline{E}_k] \nonumber \\
    \leq&\ -\gamma_{2,k}\alpha_k\P_k[E_k | C_{1,k}] \E_k[\nabla f(x_k)^Tg_k | C_{1,k} \cap E_k] \nonumber \\
        &\ - \gamma_{1,k} \alpha_k (\E_k[\nabla f(x_k)^Tg_k | C_{1,k}] - \P_k[E_k | C_{1,k}] \E_k[\nabla f(x_k)^Tg_k | C_{1,k} \cap E_k]) \nonumber \\
       =&\ -\gamma_{1,k} \alpha_k \E_k[\nabla f(x_k)^Tg_k | C_{1,k}] + (\gamma_{1,k} - \gamma_{2,k}) \alpha_k \P_k[E_k | C_{1,k}] \E_k[\nabla f(x_k)^Tg_k | C_{1,k} \cap E_k] \label{eq.case1_prod}
  \end{align}
  along with the fact that
  \bequation\label{eq.case1_noise}
    \E_k[\|s_k\|^2 | C_{1,k}] = \gamma_{1,k}^2 \alpha_k^2 \E_k[\|g_k\|^2 | C_{1,k}].
  \eequation
  In the event $C_{2,k}$, in which $\|g_k\|^{-1} \leq \gamma_{1,k}$ and $\|g_k\|^{-1} \geq \gamma_{2,k}$, one finds that
  \begin{align}
       &\ \E_k[\nabla f(x_k)^Ts_k | C_{2,k}] \nonumber \\
      =&\ -\alpha_k \E_k\left[\frac{\nabla f(x_k)^Tg_k}{\|g_k\|} \bigg| C_{2,k}\right] \nonumber \\
      =&\ -\alpha_k \P_k[E_k | C_{2,k}] \E_k\left[\frac{\nabla f(x_k)^Tg_k}{\|g_k\|} \bigg| C_{2,k} \cap E_k\right] - \alpha_k \P_k[\overline{E}_k | C_{2,k}] \E_k\left[\frac{\nabla f(x_k)^Tg_k}{\|g_k\|} \bigg| C_{2,k} \cap \overline{E}_k\right] \nonumber \\
      \leq&\ -\gamma_{2,k} \alpha_k \P_k[E_k | C_{2,k}] \E_k[\nabla f(x_k)^Tg_k | C_{2,k} \cap E_k] - \gamma_{1,k} \alpha_k \P_k[\overline{E}_k | C_{2,k}] \E_k[\nabla f(x_k)^Tg_k | C_{2,k} \cap \overline{E}_k] \nonumber \\
         =&\ - \gamma_{2,k} \alpha_k \P_k[E_k | C_{2,k}] \E_k[\nabla f(x_k)^Tg_k | C_{2,k} \cap E_k] \nonumber \\
          &\ - \gamma_{1,k} \alpha_k (\E_k[\nabla f(x_k)^Tg_k | C_{2,k}] - \P_k[E_k | C_{2,k}] \E_k[\nabla f(x_k)^Tg_k | C_{2,k} \cap E_k]) \nonumber \\
         =&\ -\gamma_{1,k}\alpha_k\E_k[\nabla f(x_k)^Tg_k | C_{2,k}] + (\gamma_{1,k} - \gamma_{2,k}) \alpha_k \P_k[E_k | C_{2,k}] \E_k[\nabla f(x_k)^Tg_k | C_{2,k} \cap E_k] \label{eq.case2_prod}
  \end{align}
  along with the fact that
  \bequation\label{eq.case2_noise}
    \E_k[\|s_k\|^2 | C_{2,k}] = \alpha_k^2 \leq \gamma_{1,k}^2 \alpha_k^2 \E_k[\|g_k\|^2 | C_{2,k}].
  \eequation
  In the event $C_{3,k}$, the algorithm yields $s_k = -\gamma_{2,k} \alpha_k g_k$, from which it follows that
  \begin{align}
        &\ \E_k[\nabla f(x_k)^Ts_k | C_{3,k}] \nonumber \\
       =&\ -\gamma_{2,k}\alpha_k\E_k[\nabla f(x_k)^Tg_k | C_{3,k}] \nonumber \\
    \leq&\ -\gamma_{2,k}\alpha_k\P_k[E_k | C_{3,k}] \E_k[\nabla f(x_k)^Tg_k | C_{3,k} \cap E_k] - \gamma_{1,k} \alpha_k \P_k[\overline{E}_k | C_{3,k}] \E_k[\nabla f(x_k)^Tg_k | C_{3,k} \cap \overline{E}_k] \nonumber \\
       =&\ -\gamma_{2,k}\alpha_k\P_k[E_k | C_{3,k}] \E_k[\nabla f(x_k)^Tg_k | C_{3,k} \cap E_k] \nonumber \\
        &\ - \gamma_{1,k} \alpha_k (\E_k[\nabla f(x_k)^Tg_k | C_{3,k}] - \P_k[E_k | C_{3,k}] \E_k[\nabla f(x_k)^Tg_k | C_{3,k} \cap E_k]) \nonumber \\
       =&\ -\gamma_{1,k} \alpha_k \E_k[\nabla f(x_k)^Tg_k | C_{3,k}] + (\gamma_{1,k} - \gamma_{2,k}) \alpha_k \P_k[E_k | C_{3,k}] \E_k[\nabla f(x_k)^Tg_k | C_{3,k} \cap E_k] \label{eq.case3_prod}
  \end{align}
  along with the fact that
  \bequation\label{eq.case3_noise}
    \E_k[\|s_k\|^2 | C_{3,k}] = \gamma_{2,k}^2 \alpha_k^2 \E_k[\|g_k\|^2 | C_{3,k}] \leq \gamma_{1,k}^2 \alpha_k^2 \E_k[\|g_k\|^2 | C_{3,k}].
  \eequation
  Combining \eqref{eq.trick1}--\eqref{eq.case3_noise}, it follows that
  \bequationn
    \baligned
          &\ \E_k[f(x_{k+1})] - f(x_k) \\
      \leq&\ - \gamma_{1,k} \alpha_k \|\nabla f(x_k)\|^2 + (\gamma_{1,k} - \gamma_{2,k}) \alpha_k \P_k[E_k] \E_k[\nabla f(x_k)^Tg_k | E_k] + \thalf \gamma_{1,k}^2 L \alpha_k^2 \E_k[\|g_k\|^2].
    \ealigned
  \eequationn
  Applying the upper bound for the last term in \eqref{eq.second_moment} and rearranging terms yields the result. 
\endproof

For some (but not all) of our convergence guarantees, we also make the following assumption.

\bassumption\label{ass.PL}
  At any $x \in \R{n}$, the Polyak-\L{}ojasiewicz condition holds with $c \in (0,\infty)$, i.e.,
  \bequation\label{eq.PL}
  2c(f(x) - f_*) \leq \|\nabla f(x)\|^2\ \ \text{for all}\ \ x \in \R{n}.
  \eequation
\eassumption

Assumptions~\ref{ass.f} and \ref{ass.PL} do not ensure that a stationary point for $f$ exists, though, when combined, they do guarantee that any stationary point for $f$ is a global minimizer of $f$.  Assumption~\ref{ass.PL} holds when $f$ is $c$-strongly convex, but it is also satisfied for other functions that are not convex.  We direct the interested reader to \cite{karimi2016linear} for a discussion on the relationship between the Polyak-\L{}ojasiewicz condition and the related \emph{error bounds}, \emph{essential strong convexity}, \emph{weak strong convexity}, \emph{restricted secant inequality}, and \emph{quadratic growth} conditions.  In short, when $f$ has a Lipschitz continuous gradient, the Polyak-\L{}ojasiewicz is the weakest of these except for the quadratic growth condition, though these two are equivalent when $f$ is convex.

We now proceed to prove convergence guarantees for \ref{alg.trish} in various cases depending on whether or not the Polyak-\L{}ojasiewicz condition (hereafter referred to as the P-L condition) holds and based on different sets of properties of the sequence of stepsizes and stochastic gradient estimates.  Our analysis covers various types of convex and nonconvex objective functions.

\subsection{P-L Condition and Constant Parameters}\label{subsec.first}

Let us first prove a convergence result for \ref{alg.trish} when the P-L condition holds and each sequence $\{\alpha_k\}$, $\{\gamma_{1,k}\}$, and $\{\gamma_{2,k}\}$ is constant.  This result appears in this section as~Theorem~\ref{th.bounded_variance_fixed_stepsize}.

Our first requirement toward proving Theorem~\ref{th.bounded_variance_fixed_stepsize} is the following lemma.

\blemma\label{lem.fg_upperbound}
  Under Assumption~\ref{ass.f_limits}, it follows that, for all $k \in \N{}$,
  \bequation\label{eq.ce_fg_upperbound}
    \P_k[E_k]\E_k[\nabla f(x_k)^T g_k | E_k] \leq h_1 + h_2\|\nabla f(x_k)\|^2
  \eequation
  for any $(h_1,h_2) \in (0,\infty) \times (0,\infty)$ such that $h_1 \geq \thalf\sqrt{M_1}$ and $h_2 \geq \thalf\sqrt{M_1} + \sqrt{M_2}$.
\elemma
\proof{Proof.}
  One finds with the law of total probability that
  \bequationn
    \baligned
      \P_k[E_k]\E_k[\nabla f(x_k)^T g_k | E_k]
        &\leq \P_k[E_k]\E_k[\|\nabla f(x_k)\| \|g_k\| | E_k] \\
        &=    \|\nabla f(x_k)\| (\P_k[E_k] \E_k[\|g_k\| | E_k]) \\
        &=    \|\nabla f(x_k)\| (\E_k[\|g_k\|] - \P_k[\overline{E}_k] \E_k[\|g_k\| | \overline{E}_k]) \\
        &\leq \|\nabla f(x_k)\| \E_k[ \|g_k\|].
    \ealigned
  \eequationn
  Then, by Jensen's Inequality, concavity of the square root, and Assumption~\ref{ass.f_limits}, one finds that
  \bequationn
    \E_k[\|g_k\|] \leq \sqrt{\E_k[\|g_k\|^2]} \leq \sqrt{M_1 + M_2 \|\nabla f(x_k) \|^2} \leq \sqrt{M_1}+ \sqrt{M_2} \|\nabla f(x_k) \|.
  \eequationn
  Therefore, by combining the inequalities above, one finds that
  \bequationn
    \baligned
      \P_k[E_k]\E_k[\nabla f(x_k)^Tg_k | E_k]
        &\leq \|\nabla f(x_k)\| (\sqrt{M_1} + \sqrt{M_2} \|\nabla f(x_k) \|) \\
        &=    \sqrt{M_1} \|\nabla f(x_k)\| + \sqrt{M_2}\|\nabla f(x_k) \|^2 \\
        &\leq \half \sqrt{M_1} (1 + \|\nabla f(x_k)\|^2) + \sqrt{M_2} \|\nabla f(x_k)\|^2 \\
        &=    \half\sqrt{M_1} + \(\half\sqrt{M_1} + \sqrt{M_2}\) \|\nabla f(x_k)\|^2,
    \ealigned
  \eequationn
  where the second inequality follows by the fact that $a \leq \thalf(1+a^2)$ for any $a \in \R{}$. 
\endproof

While the upper bound on $\E_k[\|g_k\|^2]$ stated as~\eqref{eq.second_moment} in Assumption~\ref{lem.fg_upperbound} is standard in the literature, the quantity on the left-hand side of~\eqref{eq.ce_fg_upperbound}---which Lemma~\ref{lem.fg_upperbound} shows is bounded in a similar manner---is uniquely important for our analysis.  For this reason, we feel that it is useful to provide specific examples illustrating how this quantity is bounded.  We state two related examples next.

\bexample\label{example.fixed_step}
  \textit{
  Suppose $f : \R{} \to \R{}$ and $x_k$ are given such that $\nabla f(x_k) = \mu_k \in \R{}$, where without loss of generality one can assume that $\mu_k \geq 0$.  In addition, suppose that $g_k$ follows a normal distribution with mean $\mu_k$ and variance $\sigma_k^2$.  Then,
  \begin{align*}
    \P_k[E_k]\E_k[\nabla f(x_k)^Tg_k | E_k]
      &= \mu_k \int_0^{\infty} g \tfrac{1}{\sqrt{2\pi} \sigma_k} e^{\frac{-(g-\mu_k)^2}{2\sigma_k^2}} dg \\
      &= \mu_k \int_0^{\mu_k} g \tfrac{1}{\sqrt{2\pi} \sigma_k} e^{\frac{-(g-\mu_k)^2}{2\sigma_k^2}} dg + \mu_k \int_{\mu_k}^{\infty} g \tfrac{1}{\sqrt{2\pi} \sigma_k} e^{\frac{-(g-\mu_k)^2}{2\sigma_k^2}} dg.
  \end{align*}
  Let us separately investigate these two terms on the right-hand side.  First, one finds that
  \bequationn
    \mu_k \int_0^{\mu_k} g \tfrac{1}{\sqrt{2\pi} \sigma_k} e^{\frac{-(g-\mu_k)^2}{2\sigma_k^2}} dg \leq \mu_k^2 \int_0^{\mu_k} \tfrac{1}{\sqrt{2\pi} \sigma_k} e^{\frac{-(g-\mu_k)^2}{2\sigma_k^2}} dg \leq \mu_k^2 \int_{-\infty}^{\mu_k} \tfrac{1}{\sqrt{2\pi} \sigma_k} e^{\frac{-(g-\mu_k)^2}{2\sigma_k ^2}} dg = \thalf \mu_k^2.
  \eequationn
  Second, one finds that
  \bequationn
    \baligned
      \mu_k \int_{\mu_k}^{\infty} g \tfrac{1}{\sqrt{2\pi} \sigma_k} e^{\frac{-(g-\mu_k)^2}{2\sigma_k^2}} dg
        &= \mu_k \int_0^{\infty} (t + \mu_k)\tfrac{1}{\sqrt{2\pi} \sigma_k} e^{\frac{-t^2}{2\sigma_k^2}} dt \\
        &= \mu_k \int_0^{\infty} t \tfrac{1}{\sqrt{2\pi} \sigma_k} e^{\frac{-t^2}{2\sigma_k ^2}} dt + \mu_k^2 \int_0^{\infty} \tfrac{1}{\sqrt{2\pi} \sigma_k} e^{\frac{-t^2}{2\sigma_k ^2}} dt = \mu_k \tfrac{\sigma_k}{\sqrt{2\pi}} + \thalf \mu_k^2.
    \ealigned
  \eequationn
  Thus, combining the bounds above, one finds that
  \bequationn
    \P_k[E_k]\E_k[\nabla f(x_k)^Tg_k | E_k] \leq \mu_k \frac{\sigma_k}{\sqrt{2\pi}} + \mu_k^2 \leq \(\frac{\mu_k^2 + 1}{2}\)\frac{\sigma_k}{\sqrt{2\pi}} + \mu_k^2 = \frac{\sigma_k}{2\sqrt{2\pi}} + \(1 + \frac{\sigma_k}{2\sqrt{2\pi}}\)\mu_k^2.
  \eequationn
  Overall, if $\sigma_k \leq \sigma$ for some positive $\sigma \in \R{}$ for all $k \in \N{}$, then \eqref{eq.ce_fg_upperbound} holds with
  \bequation\label{eq.example1}
    h_1 = \frac{\sigma}{2\sqrt{2\pi}}\ \ \text{and}\ \ h_2 = 1 + \frac{\sigma}{2\sqrt{2\pi}}.
  \eequation
  }
\eexample

\bexample\label{example.fixed_step_2}
  \textit{
  Suppose $f : \R{n} \to \R{}$ and $x_k$ are given such that $\nabla f(x_k) = \mu_k \in \R{n}$.  In addition, suppose that $g_k$ follows a normal distribution with mean $\mu_k$ and covariance matrix $\Sigma_k$.  Then, by Theorem~3.3.3 in \cite{tong2012multivariate}, the inner product $\nabla f(x_k)^Tg_k$ follows a normal distribution with mean $\|\mu_k\|^2$ and variance $\mu_k^T\Sigma_k\mu_k$.  Hence, following the analysis in Example~\ref{example.fixed_step}, if $\sqrt{\mu_k^T\Sigma_k\mu_k} \leq \sigma$ for some positive $\sigma \in \R{}$ for all $k \in \N{}$, then \eqref{eq.ce_fg_upperbound} holds with $h_1$ and $h_2$ from \eqref{eq.example1}.
  }
\eexample

We now prove our first theorem on the behavior of \ref{alg.trish}.

\btheorem\label{th.bounded_variance_fixed_stepsize}
  Under Assumptions~\ref{ass.f}, \ref{ass.f_limits}, and \ref{ass.PL}, and with a pair $(h_1,h_2)$ satisfying the inequalities in Lemma~\ref{lem.fg_upperbound}, suppose that \ref{alg.trish} is run with $(\gamma_{1,k},\gamma_{2,k}) = (\gamma_1,\gamma_2)$ for all $k \in \N{}$ such that $\frac{\gamma_1}{\gamma_2} < \frac{h_2}{h_2-1}$ $($meaning $\gamma_1 - h_2(\gamma_1 - \gamma_2) > 0$$)$ and with $\alpha_k = \alpha$ for all $k \in \N{}$ such that
  \bequation\label{eq.fixed_alpha_choice}
    0 < \alpha \leq \min\left\{\frac{1}{2c\theta_1}, \frac{\gamma_1 - h_2(\gamma_1 - \gamma_2)}{\gamma_1 L M_2}\right\},
  \eequation
  where
  \bequation\label{eq.theta1}
    \theta_1 = \thalf(\gamma_1 - h_2(\gamma_1 - \gamma_2)) > 0.
  \eequation
  Then, for all $k \in \N{}$, the expected optimality gap satisfies
  \bequation\label{eq.fixed_stepsize_result}
    \E[f(x_{k+1})] - f_* \leq \frac{\theta_2}{2c \alpha\theta_1} + (1 - 2c \alpha \theta_1)^{k-1} \(f(x_1) - f_* - \frac{\theta_2}{2c \alpha \theta_1}\) \xrightarrow{k\rightarrow\infty} \frac{\theta_2}{2c \alpha \theta_1},
  \eequation
  where
  \bequation\label{eq.theta2}
    \theta_2 = h_1(\gamma_1 - \gamma_2)\alpha + \thalf \gamma_1^2 L M_1 \alpha^2 > 0. 
  \eequation
\etheorem
\proof{Proof.}
  Combining the results of Lemmas~\ref{lem.basic} and \ref{lem.fg_upperbound}, it follows that, for all $k \in \N{}$,
  \bequation\label{eq.fixed_stepsize_bound}
    \baligned
      \E_k[f(x_{k+1})] - f(x_k)
        \leq&\ -\gamma_1 \alpha (1 - \thalf \gamma_1 L M_2 \alpha) \|\nabla f(x_k)\|^2 \\
            &\ + (\gamma_1 - \gamma_2) \alpha (h_1 + h_2\|\nabla f(x_k)\|^2) + \thalf \gamma_1^2 L M_1 \alpha^2.
    \ealigned
  \eequation
  Therefore, with $(\theta_1,\theta_2)$ defined in \eqref{eq.theta1}/\eqref{eq.theta2}, it follows with \eqref{eq.PL} that, for all $k \in \N{}$,
  \bequationn
    \baligned
      \E_k[f(x_{k+1})] - f(x_k)
        &\leq -\alpha \theta_1 \|\nabla f(x_k)\|^2 + \theta_2 \\
        &\leq - 2c \alpha \theta_1 (f(x_k) - f_*) + \theta_2.
    \ealigned
  \eequationn
  Adding and subtracting $f_*$ on the left-hand side, taking total expectations, and rearranging yields
  \bequationn
    \baligned
      \E[f(x_{k+1})] - f_*
        &\leq (1 - 2c \alpha \theta_1) (\E[f(x_{k})] - f_*) + \theta_2 \\
        &=    \frac{\theta_2}{2c\alpha\theta_1} + (1 - 2c \alpha \theta_1) (\E[f(x_{k})] - f_*) + \theta_2 - \frac{\theta_2}{2c\alpha\theta_1} \\
        &=    \frac{\theta_2}{2c\alpha\theta_1} + (1 - 2c \alpha \theta_1) \(\E[f(x_{k})] - f_* - \frac{\theta_2}{2c\alpha\theta_1}\).
    \ealigned
  \eequationn
  Since $1 - 2c\alpha\theta_1 \in (0,1)$, this represents a contraction inequality.  Applying the result repeatedly through iteration $k \in \N{}$, one obtains the desired result.
\endproof

It is worthwhile to compare the result of Theorem~\ref{th.bounded_variance_fixed_stepsize} with a corresponding result known to hold for a straightforward SG method.  For example, from \cite[Thm.~4.6]{BottCurtNoce18} with our notation, it is known that for an SG method with fixed stepsize $\alpha = \frac{1}{LM_2}$ an upper bound for the expected optimality gap converges to $\frac{\alpha L M_1}{2c} = \frac{M_1}{2cM_2}$.  On the other hand, the analysis in Theorem~\ref{th.bounded_variance_fixed_stepsize} shows that \ref{alg.trish} with $\alpha = \frac{\gamma_1-h_2(\gamma_1-\gamma_2)}{\gamma_1LM_2}$ (which may occur, e.g., if $c \approx 0$) yields an upper bound for the expected optimality gap that converges to
\bequation\label{eq.discussion}
  \frac{h_1 (\gamma_1 - \gamma_2) + \thalf \gamma_1^2 L M_1 \alpha}{c(\gamma_1 - h_2(\gamma_1 - \gamma_2))} = \frac{h_1 (\gamma_1 - \gamma_2)}{c(\gamma_1 - h_2(\gamma_1 - \gamma_2))} + \frac{\gamma_1 M_1}{2 c M_2}.
\eequation
We can now make a couple of observations.  On one hand, if $h_1 \approx \thalf \sqrt{M_1}$ and $h_2 \approx M_2 \approx 1$, then the condition that $\tfrac{\gamma_1}{\gamma_2} < \tfrac{h_2}{h_2-1}$ essentially does not restrict $(\gamma_1,\gamma_2)$, in which case \eqref{eq.discussion} is approximately
\bequationn
  \frac{\sqrt{M_1}(\gamma_1 - \gamma_2)}{2c\gamma_2} + \frac{\gamma_1 M_1}{2c}.
\eequationn
This quantity is less than $\frac{M_1}{2c}$, i.e., the approximate bound for SG, if, e.g., the parameters satisfy $\gamma_1 \in (0,1)$ with $\gamma_2 \geq \tfrac{\gamma_1}{1 + (1-\gamma_1)\sqrt{M_1}} \in (0,\gamma_1)$.  On the other hand, if $h_1 \approx \thalf \sqrt{M_1}$ and $h_2 \approx \thalf \sqrt{M_1} + \sqrt{M_2}$ with $M_1 \gg 0$, then the condition that $\tfrac{\gamma_1}{\gamma_2} < \tfrac{h_2}{h_2-1}$ essentially requires that $\gamma_1 \approx \gamma_2$, in which case the bound \eqref{eq.discussion} is approximately~$\tfrac{\gamma_1 M_1}{2cM_2}$, which is less than the bound for SG if $\gamma_1 \in (0,1)$.  Overall, while we are not necessarily recommending that one employes \ref{alg.trish} with the parameter settings mentioned in this discussion, we have at least been able to demonstrate in both of these cases that \ref{alg.trish} can possess an asymptotic bound on the expected optimality gap that is on par with that for SG.  (For a detailed discussion on how to choose $(\gamma_1,\gamma_2)$ in practice, see~\S\ref{subsec.parameter}.)

Besides the conclusions of the previous paragraph, the result of Theorem~\ref{th.bounded_variance_fixed_stepsize} points to fundamental differences between \ref{alg.trish} and SG for certain choices of the input parameters.  In particular, the result in \cite[Thm.~4.6]{BottCurtNoce18} points to a well-known trade-off for SG with a fixed stepsize: If a relatively large stepsize is employed, then the rate to achieve the asymptotic expected optimality gap involves a better constant at the sake of the upper bound on the gap being relatively large, i.e., $\tfrac{\alpha L M_1}{2c}$, which is proportional to the stepsize $\alpha$.  On the other hand, one can achieve a smaller upper bound on the expected optimality gap with a smaller $\alpha$, but at the cost of a worse constant in the rate to achieve that gap.  A similar conclusion can be derived from \eqref{eq.fixed_stepsize_result} for \ref{alg.trish}: One can control the constant $(1 - 2c\alpha\theta_1)$ by the choice of $\alpha$.  However, the effect of $\alpha$ on the expected optimality gap is not exactly the same for \ref{alg.trish} as for an SG method.  This can be seen in the fact that the left-hand side of~\eqref{eq.discussion} involves one term that decreases with $\alpha$, but another term that does not.  That said, one can compensate for this in \ref{alg.trish} if one ties the difference $\gamma_1 - \gamma_2$ to the stepsize $\alpha$.  This idea can be seen in the first of our two theorems in the next subsection.

\subsection{P-L Condition and Sublinearly Diminishing Stepsizes}

Let us now consider the behavior of \ref{alg.trish} when the P-L condition holds and diminishing stepsizes are employed.  Our first theorem in this setting, which makes the same assumptions as Theorem~\ref{th.bounded_variance_fixed_stepsize}, but involves different parameter choices, is the following.  (The parameter choices in the theorem could be generalized even further.  However, we have made certain choices---e.g., to have $\{\gamma_1\}$ be constant---for some amount of simplicity in the proof while still maintaining generality.  One could prove a similar result with $\{\gamma_2\}$ constant instead, or with neither $\{\gamma_1\}$ nor $\{\gamma_2\}$ constant, as long as the sequence $\{\gamma_{1,k} - \gamma_{2,k}\}$ is proportional to $\alpha_k$, as it is in the following theorem.)

\begin{theorem}\label{th.diminishing_stepsize}
  Under Assumptions~\ref{ass.f}, \ref{ass.f_limits}, and \ref{ass.PL}, and with a pair $(h_1,h_2)$ satisfying the inequalities in Lemma~\ref{lem.fg_upperbound}, suppose that \ref{alg.trish} is run with $\gamma_{1,k} = \gamma_1 > 0$, $\gamma_{2,k} = \gamma_1(1 - \thalf\eta\alpha_k)$ for $\eta \in (0,1)$, and
  \bequation\label{eq.alpha_dim}
    \alpha_k = \frac{a}{b+k}\ \text{for some}\ a \in \(\frac{1}{c\gamma_1},\frac{b+1}{c\gamma_1}\)\ \text{and}\ b > 0\ \text{with}\ \alpha_1 \in \left(0,\min\left\{\frac{1}{\eta},\frac{1}{\eta h_2 + \gamma_1 L M_2}\right\}\right]
  \eequation
  for all $k \in \N{}$.  Then, for all $k \in \N{}$, the expected optimality gap satisfies
  \bequation\label{eq.sublinear_1}
    \E[f(x_k)] - f_* \leq \frac{\phi}{b+k},  
  \eequation
  where
  \begin{align}
    \phi &= \max\left\{ \frac{a^2 \delta}{ac\gamma_1-1}, (b+1) (f(x_1)-f_*) \right\} > 0 \label{eq.phi} \\
    \text{and}\ \ \delta &= \thalf \gamma_1(\eta h_1 + \gamma_1 L M_1) > 0. \label{eq.delta}
  \end{align}
\end{theorem}
\proof{Proof.}
  First observe that the restrictions on $\{\alpha_k\}$ in \eqref{eq.alpha_dim} ensure that $\gamma_{2,k} > 0$, $\gamma_1 - \gamma_{2,k} = \thalf\gamma_1\eta\alpha_k$, and $1 - \thalf(\eta h_2 + \gamma_1LM_2)\alpha_k \geq \thalf$ for all $k \in \N{}$.  Thus, similar to the proof of Theorem \ref{th.bounded_variance_fixed_stepsize}, for all $k \in \N{}$,
  \bequationn
    \baligned
      \E_k[f(x_{k+1})] - f(x_k)
        \leq&\ - \gamma_1 \alpha_k(1 - \thalf \gamma_1 L M_2 \alpha_k) \|\nabla f(x_k)\|^2 \\
            &\ + (\gamma_1 - \gamma_{2,k}) \alpha_k (h_1 + h_2\|\nabla f(x_k)\|^2) + \thalf \gamma_1^2 L M_1 \alpha_k^2 \\
           =&\ - \gamma_1 \alpha_k(1 - \thalf \gamma_1 L M_2 \alpha_k) \|\nabla f(x_k)\|^2 \\
            &\ + \thalf \gamma_1 \eta \alpha_k^2 (h_1 + h_2\|\nabla f(x_k)\|^2) + \thalf \gamma_1^2 L M_1 \alpha_k^2 \\
           =&\ - \gamma_1 \alpha_k(1 - \thalf (\eta h_2 + \gamma_1 L M_2) \alpha_k) \|\nabla f(x_k)\|^2 + \thalf \gamma_1 (\eta h_1 + \gamma_1 L M_1) \alpha_k^2 \\
        \leq&\ - \thalf \gamma_1 \alpha_k \|\nabla f(x_k)\|^2 + \thalf \gamma_1 (\eta h_1 + \gamma_1 L M_1) \alpha_k^2.
    \ealigned
  \eequationn
  Therefore, with $\delta$ defined in \eqref{eq.delta}, it follows with \eqref{eq.PL} that, for all $k \in \N{}$,
  \begin{align}
    \E_k[f(x_{k+1})] - f(x_k)
      &\leq -\thalf \gamma_1 \alpha_k \|\nabla f(x_k)\|^2 + \delta \alpha_k^2 \label{eq.for_later_1} \\
      &\leq -c \gamma_1 \alpha_k (f(x_k) - f_*) + \delta \alpha_k^2. \nonumber
  \end{align}
  Adding and subtracting $f_*$ on the left-hand side, taking total expectations, and rearranging yields
  \bequation\label{eq.jjj}
    \E[f(x_{k+1})] - f_* \leq (1 - c\gamma_1 \alpha_k) (\E[f(x_{k})] - f_*) + \delta \alpha_k^2.
  \eequation
  
  Let us now prove \eqref{eq.sublinear_1} by induction.  First, for $k=1$, the inequality holds by the definition of $\phi$.  Now suppose that \eqref{eq.sublinear_1} holds up to $k \in \N{}$; then, for $k+1$, one finds from above that
  \begin{align*}
    \E[f(x_{k+1})] - f_*
      &\leq (1 - c\gamma_1 \alpha_k) (\E[f(x_k)] - f_*) + \delta \alpha_k^2 \\
      &=    \(1 - \frac{ac \gamma_1}{b+k}\)(\E[f(x_{k})] - f_*) + \frac{a^2\delta}{(b+k)^2}\\
      &\leq \(1 - \frac{ac \gamma_1}{b+k}\) \frac{\phi}{b+k} + \frac{a^2\delta}{(b+k)^2} \\
      &=    \frac{(b+k)\phi}{(b+k)^2} - \frac{ac\gamma_1 \phi}{(b+k)^2} + \frac{a^2\delta}{(b+k)^2}\\
      &=    \frac{(b+k-1)\phi}{(b+k)^2} - \frac{(ac\gamma_1-1) \phi}{(b+k)^2} + \frac{a^2\delta}{(b+k)^2}\\
      &\leq \frac{(b+k-1)\phi}{(b+k)^2} \leq \frac{\phi}{b+k+1},
  \end{align*}
  where the last two inequalities follow from the definition of $\phi$ and since $(b+k-1)(b+k+1) \leq (b+k)^2$, respectively.  The desired conclusion now follows from this inductive argument. 
\endproof

As one might predict from the discussion at the end of \S\ref{subsec.first}, in Theorem~\ref{th.diminishing_stepsize} we have been able to prove sublinear convergence of the expected optimality gap by tying the rate that $\{\gamma_{1,k} - \gamma_{2,k}\}$ vanishes to the rate that $\{\alpha_k\}$ vanishes; in particular, both the differences and the stepsizes diminish sublinearly, as is the case in similar results for SG methods.

One might also be interested in the behavior of \ref{alg.trish} when the sequences $\{\gamma_{1,k}\}$ and $\{\gamma_{2,k}\}$ are constant while only the stepsizes decrease sublinearly.  For example, this might be of interest since otherwise there are additional parameters to estimate and/or to tune.  In the remainder of this subsection, we prove a sublinear convergence result under this setting.  However, achieving sublinear convergence in this setting requires the following assumption, which can be viewed as a strengthening of \eqref{eq.second_moment} from Assumption~\ref{ass.f_limits}.

\bassumption\label{ass.fg_upperbound2}
  There exists a pair $(M_3,M_4) \in (0,\infty) \times (0,\infty)$ (independent of $k$) such that, for all $k \in \N{}$, the squared norm of $g_k$ satisfies
  \bequation\label{eq.second_moment2}
    \E_k[\|g_k\|^2] \leq M_3 \alpha_k^2 + M_4 \|\nabla f(x_k) \|^2.
  \eequation
\eassumption

One finds that Assumption~\ref{ass.fg_upperbound2} can be satisfied under reasonable conditions in practice if one employs mini-batch stochastic gradient estimates with sample sizes that increase with $k$; see, e.g., \cite{FrieSchm12}.  For example, in the context of problem~\eqref{prob.f}, suppose that
\bequation\label{eq.mini-batch-sg}
  g_k = \frac{1}{|\Scal_k|} \sum_{j\in\Scal_k} \nabla_x F(x_k,\xi_{k,j}),
\eequation
where the values $\{\xi_{k,j}\}_{j\in\Scal_k}$ are drawn independently according to the distribution of $\xi$.  If one assumes that the variance of each $\nabla_x F(x_k,\xi_{k,j})$ is equal and bounded by $M \in (0,\infty)$, then for arbitrary $j \in \Scal_k$ it follows (see, e.g., \cite{Freu62}) that
\bequation\label{eq.mini-batch-variance}
  \E_k[\|g_k\|^2] - \|\nabla f(x_k)\|^2 \leq \frac{M}{|\Scal_k|}.
\eequation
Hence, \eqref{eq.second_moment2} holds with $M_3=M$ and $M_4=1$ if one chooses $|\Scal_k| = \alpha_k^{-2}$.  (In Theorem~\ref{th.diminishing_variance_diminishing_stepsize} below, the result requires $\alpha_k = \Theta(\tfrac1k)$, in which case one can employ $|\Scal_k| = \Theta(k^2)$.)

An important consequence of Assumption~\ref{ass.fg_upperbound2} is the following, which strengthens Lemma~\ref{lem.fg_upperbound}.

\blemma\label{lem.fg_upperbound2}
  Under Assumption~\ref{ass.fg_upperbound2}, it follows that, for all $k \in \N{}$,
  \bequation\label{eq.ce_fg_upperbound2}
    \P_k[E_k] \E_k[\nabla f(x_k)^T g_k | E_k] \leq h_3 \alpha_k + h_4\|\nabla f(x_k)\|^2
  \eequation
  for any $(h_3,h_4) \in (0,\infty) \times (0,\infty)$ such that $h_3 \geq \thalf \sqrt{M_3}$ and $h_4 \geq \thalf \sqrt{M_3} (\displaystyle\max_{k\in\N{}} \alpha_k) + \sqrt{M_4}$.
\elemma
\proof{Proof.}
  By Jensen's Inequality, concavity of the square root, and Assumption~\ref{ass.fg_upperbound2}, one finds that
  \bequationn
    \E_k[\|g_k\|] \leq \sqrt{\E_k[\|g_k\|^2]} \leq \sqrt{M_3 \alpha_k^2 + M_4 \|\nabla f(x_k) \|^2} \leq \sqrt{M_3} \alpha_k+ \sqrt{M_4} \|\nabla f(x_k)\|.
  \eequationn
  The result then follows using the same line of argument as used in the proof of Lemma~\ref{lem.fg_upperbound}. 
\endproof

The following examples parallel Examples~\ref{example.fixed_step} and \ref{example.fixed_step_2}, but illustrate the attainment of~\eqref{eq.ce_fg_upperbound2}.

\bexample
  \textit{
  Consider the scenario in Example~\ref{example.fixed_step}.  Then, if $\sigma_k \leq \alpha_k$ for all $k \in \N{}$ with $\alpha_k \leq \alpha$ for some $\alpha \in (0,\infty)$ for all $k \in \N{}$, it follows that \eqref{eq.ce_fg_upperbound2} holds with
  \bequation\label{eq.example2}
    h_3 = \frac{1}{2\sqrt{2\pi}}\ \ \text{and}\ \ h_4 = 1 + \frac{\alpha}{2\sqrt{2\pi}}.
  \eequation
  }
\eexample

\bexample
  \textit{
  Consider the scenario in Example~\ref{example.fixed_step_2}.  Then, if $\sqrt{\mu_k^T\Sigma_k\mu_k} \leq \alpha_k$ for all $k \in \N{}$ with $\alpha_k \leq \alpha$ for some $\alpha \in (0,\infty)$ for all $k \in \N{}$, it follows that \eqref{eq.ce_fg_upperbound2} holds with $h_3$ and $h_4$ from~\eqref{eq.example2}.
  }
\eexample

Our next theorem on the behavior of \ref{alg.trish} is now proved as the following.  (For the result, we include Assumptions~\ref{ass.f_limits} and \ref{ass.fg_upperbound2} for convenience since, in our proof, we employ results that we have proved using each of these assumptions.  Notice, however, that the bound~\eqref{eq.second_moment} in Assumption~\ref{ass.f_limits} holds under Assumption~\ref{ass.fg_upperbound2} if one considers $M_1 \geq M_3(\displaystyle\max_{k\in\N{}} \alpha_k^2)$ and $M_2 = M_4$.)

\begin{theorem}\label{th.diminishing_variance_diminishing_stepsize}
  Under Assumptions~\ref{ass.f}, \ref{ass.f_limits}, \ref{ass.PL}, and \ref{ass.fg_upperbound2}, and with a pair $(h_3,h_4)$ satisfying the inequalities in Lemma~\ref{lem.fg_upperbound2}, suppose that \ref{alg.trish} is run with $\gamma_1 > \gamma_2 > 0$ such that $\frac{\gamma_1}{\gamma_2} < \frac{h_4}{h_4-1}$ $($meaning $\gamma_1 - h_4(\gamma_1 - \gamma_2) > 0$$)$, and with, for all $k \in \N{}$,
  \begin{equation*}
    \alpha_k = \frac{a}{b+k}\ \ \text{for some}\ \ a \in \(\frac{1}{2c\beta_1},\frac{b+1}{2c\beta_1}\)\ \ \text{and}\ \ b > 0\ \ \text{such that}\ \ \alpha_1 \in \left(0,\frac{\gamma_1 - h_4(\gamma_1 - \gamma_2)}{\gamma_1 L M_2}\right],
  \end{equation*}
  where  
  \bequation\label{eq.beta1}
    \beta_1 = \thalf(\gamma_1 - h_4(\gamma_1 -\gamma_2)) > 0.
  \eequation
  Then, for all $k \in \N{}$, the expected optimality gap satisfies
  \bequation\label{eq.sublinear}
    \E[f(x_k)] - f_* \leq \frac{\nu}{b+k},  
  \eequation
  where
  \begin{align}
    \nu &= \max\left\{ \frac{a^2 \beta_2}{2ac\beta_1-1}, (b+1) (f(x_1)-f_*) \right\} > 0 \label{eq.nu} \\
    \text{and}\ \ \beta_2 &= h_3(\gamma_1 - \gamma_2) + \thalf \gamma_1^2 L M_1 > 0. \label{eq.beta2}
  \end{align}
\end{theorem}
\proof{Proof.}
  Similar to the proof of Theorem \ref{th.bounded_variance_fixed_stepsize}, for all $k \in \N{}$,
  \bequationn
    \baligned
      \E_k[f(x_{k+1})] - f(x_k)
        \leq&\ - \gamma_1 \alpha_k(1 - \thalf \gamma_1 L M_2 \alpha_k) \|\nabla f(x_k)\|^2 \\
            &\ + (\gamma_1 - \gamma_2) \alpha_k (h_3\alpha_k + h_4\|\nabla f(x_k)\|^2) + \thalf \gamma_1^2 L M_1 \alpha_k^2.
    \ealigned
  \eequationn
  Therefore, with $(\beta_1,\beta_2)$ defined in \eqref{eq.beta1}/\eqref{eq.beta2}, it follows with \eqref{eq.PL} that, for all $k \in \N{}$,
  \begin{align}
    \E_k[f(x_{k+1})] - f(x_k)
      &\leq -\beta_1 \alpha_k \|\nabla f(x_k)\|^2 + \beta_2 \alpha_k^2 \label{eq.for_later_2} \\
      &\leq -2c \beta_1 \alpha_k (f(x_k) - f_*) + \beta_2 \alpha_k^2. \nonumber
  \end{align}
  Adding and subtracting $f_*$ on the left-hand side, taking total expectations, and rearranging yields
  \bequationn
    \E[f(x_{k+1})] - f_* \leq (1 - 2c\beta_1 \alpha_k) (\E[f(x_{k})] - f_*) + \beta_2 \alpha_k^2.
  \eequationn
  Using this inequality, which has the same form as \eqref{eq.jjj}, one can apply the same inductive argument as in the remainder of the proof of Theorem~\ref{th.diminishing_stepsize} to achieve the desired result.
\endproof

Overall, we have proved two theorems for \ref{alg.trish} when diminishing stepsizes are employed.  If the sequence $\{\gamma_{k,1} - \gamma_{k,2}\}$ diminishes proportionally with $\{\alpha_k\}$, then sublinear convergence of the expected optimality gap is achieved under the same assumptions as needed for such a result for an SG method.  We followed this with a result for the case when $\{\gamma_{k,1} - \gamma_{k,2}\}$ is constant, in which case a sublinear convergence result for the expected optimality gap requires that the stochastic gradient estimates satisfy  Assumption~\ref{ass.fg_upperbound2}.

\subsection{P-L Condition, Constant Parameters, and Linearly Decreasing Variance}

Let us now prove a convergence result for \ref{alg.trish} when the P-L condition holds, each sequence $\{\alpha_k\}$, $\{\gamma_{1,k}\}$, and $\{\gamma_{2,k}\}$ is constant, and the stochastic gradients satisfy the following assumption.

\bassumption\label{ass.f_limits_2}
  There exist constants $(M_5,\zeta) \in (0,\infty) \times (0,1)$ such that
  \bequation\label{eq.second_moment_2}
    \E_k[\|g_k\|^2] \leq M_5 \zeta^{k-1} + \|\nabla f(x_k)\|^2.
  \eequation
\eassumption

The achievement of linear convergence of the expected optimality gap for SG also requires increasingly accurate gradient estimates along the lines required in Assumption~\ref{ass.f_limits_2}; see, e.g., \cite{BottCurtNoce18}.  One finds that Assumption~\ref{ass.f_limits_2} can be satisfied under reasonable conditions in practice if one employs mini-batch stochastic gradient estimates with sample sizes that increase with $k$.  For example, using estimates as in \eqref{eq.mini-batch-sg} and under the same conditions as led to \eqref{eq.mini-batch-variance}, one finds that \eqref{eq.second_moment_2} holds if the sample sizes increase geometrically, e.g., $|\Scal_k| = \lceil \tau^{k-1} \rceil$ for some $\tau \in (1,\infty)$.

Our main result in this section, namely, Theorem~\ref{th.diminishing_variance_fixed_stepsize}, requires the following.

\blemma\label{lem.fg_upperbound3}
  Under the Assumption~\ref{ass.f_limits_2}, it follows that, for all $k \in \N{}$,
  \bequation\label{eq.ce_fg_upperbound3}
    \P_k[E_k] \E_k[\nabla f(x_k)^T g_k | E_k] \leq h_5 \lambda^{k-1} + h_6\|\nabla f(x_k)\|^2
  \eequation
  for any $(h_5,h_6) \in (0,\infty) \times (0,\infty) \times (0,1)$ such that $h_5 \geq \thalf \sqrt{M_5}$, $1 + h_6 \geq \thalf\sqrt{M_5}$, and $\lambda \geq \sqrt{\zeta}$.
\elemma
\proof{Proof.}
  By Jensen's inequality, concavity of the square root, and Assumption~\ref{ass.f_limits_2}, one finds that
  \bequation
    \E_k[\|g_k\|] \leq  \sqrt{\E_k[ \|g_k\|^2] } \leq \sqrt{M_5 \zeta^{k-1}+ \|\nabla f(x_k) \|^2} \leq \sqrt{M_5} (\sqrt{\zeta})^{k-1} + \|\nabla f(x_k)\|.
  \eequation
  The result then follows using the same line of argument as used in the proof of Lemma~\ref{lem.fg_upperbound}. 
\endproof

The following examples parallel Examples~\ref{example.fixed_step} and \ref{example.fixed_step_2}, but illustrate the attainment of~\eqref{eq.ce_fg_upperbound3}.

\bexample
  \textit{
  Consider the scenario in Example~\ref{example.fixed_step}.  Then, since \eqref{eq.second_moment_2} implies that $\sigma_k^2 \leq M_3\zeta^{k-1}$ for all $k \in \N{}$, it follows along with the fact that $\zeta \in (0,1)$ that
  \bequationn
    \baligned
      \P_k[E_k] \E_k[\nabla f(x_k)^Tg_k | E_k]
        &\leq \frac{\sigma_k}{2\sqrt{2\pi}} + \(1 + \frac{\sigma_k}{2\sqrt{2\pi}}\)\mu_k^2 \\
        &\leq \frac{\sqrt{M_3}}{2\sqrt{2\pi}} (\sqrt{\zeta})^{k-1} + \(1 + \frac{\sqrt{M_3}}{2\sqrt{2\pi}}\) \mu_k^2.
    \ealigned
  \eequationn
  Hence, it follows that \eqref{eq.ce_fg_upperbound3} holds with
  \bequation\label{eq.example3}
    h_5 = \frac{\sqrt{M_3}}{2\sqrt{2\pi}},\ \ h_6 = 1 + \frac{\sqrt{M_3}}{2\sqrt{2\pi}},\ \ \text{and}\ \ \lambda = \sqrt{\zeta}.
  \eequation
  }
\eexample

\bexample
  \textit{
  Consider the scenario in Example~\ref{example.fixed_step_2}.  Then, with $\sqrt{\mu_k^T\Sigma_k\mu_k} \leq M_3\zeta^{k-1}$ for all $k \in \N{}$, it follows that \eqref{eq.ce_fg_upperbound3} holds with $h_5$, $h_6$, and $\lambda$ from \eqref{eq.example3}.
  }
\eexample

Our next theorem on the behavior of \ref{alg.trish} is now proved as the following.  (For the result, we include Assumptions~\ref{ass.f_limits} and \ref{ass.f_limits_2} for convenience since, in our proof, we employ results that we have proved using each of these assumptions.  Notice, however, that the bound~\eqref{eq.second_moment} in Assumption~\ref{ass.f_limits} holds under Assumption~\ref{ass.f_limits_2} if one considers $M_1 \geq M_5$, $M_2 \geq 1$, and any $\zeta \in (0,1)$.)

\btheorem\label{th.diminishing_variance_fixed_stepsize}
  Under Assumptions~\ref{ass.f},  \ref{ass.f_limits}, \ref{ass.PL}, and \ref{ass.f_limits_2}, and with a tuple $(h_5,h_6,\lambda)$ satisfying the inequalities in Lemma~\ref{lem.fg_upperbound3}, suppose that \ref{alg.trish} is run with $\gamma_1 > \gamma_2 > 0$ such that $\frac{\gamma_1}{\gamma_2} < \frac{h_6}{h_6-1}$ $($meaning $\gamma_1 - h_6(\gamma_1 - \gamma_2) > 0$$)$, and with $\alpha_k = \alpha$ for all $k \in \N{}$ such that
  \bequation\label{eq.last_alpha}
    0 < \alpha \leq \min\left\{\frac{\gamma_1 - h_6(\gamma_1 - \gamma_2)}{\gamma_1^2 L} , \frac{1}{c\kappa_1} \right\},
  \eequation
  where
  \bequation\label{eq.kappa1}
    \kappa_1 := \thalf (\gamma_1 - h_6(\gamma_1 - \gamma_2)) > 0.
  \eequation
  Then, for all $k \in \N{}$, the expected optimality gap satisfies
  \bequation\label{eq.dss.fixed_stepsize_result}
    \E[f(x_k)] - f_* \leq \omega \rho^{k-1},
  \eequation
  where
  \begin{align}
    \kappa_2 &:= h_5(\gamma_1 - \gamma_2) + \thalf \gamma_1^2 \alpha L M_3 > 0, \label{eq.kappa2} \\
    \omega   &:= \max\{f(x_1) - f_*,\tfrac{\kappa_2}{c\kappa_1}\} > 0, \nonumber \\ \text{and}\ \ 
    \rho     &:= \max\{1 - \alpha c \kappa_1 ,\lambda,\zeta\} \in (0,1). \nonumber
  \end{align}
\etheorem
\proof{Proof.}
  As in the proof of Lemma~\ref{lem.basic}, it follows with \eqref{eq.second_moment_2} and \eqref{eq.ce_fg_upperbound3} that, for all $k \in \N{}$,
  \bequationn
    \baligned
          &\ \E_k[f(x_{k+1})] - f(x_k) \\
      \leq&\ - \alpha \gamma_1 \|\nabla f(x_k)\|^2 + (\gamma_1 - \gamma_2) \alpha \P_k[E_k] \E_k[\nabla f(x_k)^Tg_k | E_k] + \thalf \gamma_1^2 L \alpha^2 \E_k[\|g_k\|^2] \\
      \leq&\ - \alpha \gamma_1 \|\nabla f(x_k)\|^2 + (\gamma_1 - \gamma_2) \alpha (h_5 \lambda^{k-1} + h_6\|\nabla f(x_k)\|^2) + \thalf \gamma_1^2 L \alpha^2 (M_3\zeta^{k-1} + \|\nabla f(x_k)\|^2) \\
         =&\ - \alpha (\gamma_1 - h_6 (\gamma_1 - \gamma_2) - \thalf \gamma_1^2 L \alpha) \|\nabla f(x_k)\|^2 + (\gamma_1 - \gamma_2) \alpha h_5 \lambda^{k-1} + \thalf \gamma_1^2 L \alpha^2 M_3 \zeta^{k-1} \\
      \leq&\ - \thalf \alpha (\gamma_1 - h_6(\gamma_1 - \gamma_2)) \|\nabla f(x_k)\|^2 + (\gamma_1 - \gamma_2) \alpha h_5 \lambda^{k-1} + \thalf \gamma_1^2 L \alpha^2 M_3 \zeta^{k-1}.
    \ealigned
  \eequationn
  Therefore, with $(\kappa_1,\kappa_2)$ defined in \eqref{eq.kappa1}/\eqref{eq.kappa2}, it follows with \eqref{eq.PL} that, for all $k \in \N{}$,
  \bequationn
    \baligned
      \E_k [f(x_{k+1})]
        &\leq f(x_k) - \alpha \kappa_1 \|\nabla f(x_k)\|^2 + \alpha\kappa_2\max\{\lambda,\zeta\}^{k-1} \\
        &\leq f(x_k) - 2\alpha c \kappa_1 (f(x_k) - f_*) + \alpha \kappa_2 \max\{\lambda,\zeta\}^{k-1},
    \ealigned
  \eequationn
  from which it follows that
  \bequationn
    \E[f(x_{k+1})] - f_* \leq (1 - 2\alpha c \kappa_1)(\E[f(x_k)] - f_*) + \alpha \kappa_2 \max\{\lambda,\zeta\}^{k-1}.
  \eequationn
  
  Let us now prove \eqref{eq.dss.fixed_stepsize_result} by induction.  First, for $k = 1$, the inequality follows by the definition of~$\omega$.  Then, assuming the inequality holds true for $k \in \N{}$, one finds that
  \bequationn
    \baligned
       \E[f(x_{k+1})] - f_* &\leq (1 - 2\alpha c \kappa_1)(\E[f(x_k)] - f_*) + \alpha \kappa_2 \max\{\lambda,\zeta\}^{k-1} \\
       &\leq (1 - 2\alpha c \kappa_1) \omega \rho^{k-1} + \alpha \kappa_2 \max\{\lambda,\zeta\}^{k-1} \\
       &=    \omega \rho^{k-1} \(1 - 2\alpha c \kappa_1 + \frac{\alpha \kappa_2}{\omega} \(\frac{\max\{\lambda,\zeta\}}{\rho}\)^{k-1}\) \\
       &\leq \omega \rho^{k-1} \(1 - 2\alpha c \kappa_1 + \frac{\alpha \kappa_2}{\omega}\) \\
       &\leq \omega \rho^{k-1} (1 - \alpha c \kappa_1) \\
       &\leq \omega \rho^k,
    \ealigned
  \eequationn
  which proves that the conclusion holds for $k + 1$, as desired. 
\endproof

\subsection{No P-L Condition and Constant Parameters}

Let us now consider the behavior of \ref{alg.trish} when the P-L condition does not hold.  Our first such result involves the use of constant $\{\gamma_{1,k}\}$, $\{\gamma_{2,k}\}$, and $\{\alpha_k\}$.

\btheorem\label{th.bounded_variance_fixed_stepsize_no_PL}
  Under Assumptions~\ref{ass.f} and \ref{ass.f_limits} and with a pair $(h_1,h_2)$ satisfying the inequalities in Lemma~\ref{lem.fg_upperbound}, suppose that \ref{alg.trish} is run with $(\gamma_{1,k},\gamma_{2,k}) = (\gamma_1,\gamma_2)$ for all $k \in \N{}$ such that $\frac{\gamma_1}{\gamma_2} < \frac{h_2}{h_2-1}$ $($meaning $\gamma_1 - h_2(\gamma_1 - \gamma_2) > 0$$)$ and with $\alpha_k = \alpha$ for all $k \in \N{}$ such that
  \bequationn
    0 < \alpha \leq \frac{\gamma_1 - h_2(\gamma_1 - \gamma_2)}{\gamma_1 L M_2}.
  \eequationn
  Then, with $(\theta_1,\theta_2)$ defined in \eqref{eq.theta1}/\eqref{eq.theta2}, it follows that, for all $K \in \N{}$,
  \bsubequations\label{eq.nonconvex}
    \begin{align}
      \E\left[ \sum_{k=1}^K \|\nabla f(x_k)\|^2 \right] \leq&\ \frac{K\theta_2}{\alpha\theta_1} + \frac{f(x_1) - f_*}{\alpha\theta_1} \label{eq.nonconvex_1} \\
      \text{and}\ \ \E\left[ \frac{1}{K} \sum_{k=1}^K \|\nabla f(x_k)\|^2 \right] \leq&\ \frac{\theta_2}{\alpha\theta_1} + \frac{f(x_1) - f_*}{K\alpha\theta_2} \xrightarrow{K\to\infty} \frac{\theta_2}{\alpha\theta_1}. \label{eq.nonconvex_2}
    \end{align}
  \esubequations
\etheorem
\proof{Proof.}
  As in the proof of Theorem~\ref{th.bounded_variance_fixed_stepsize}, combining the results of Lemmas~\ref{lem.basic} and \ref{lem.fg_upperbound}, it follows that the inequality~\eqref{eq.fixed_stepsize_bound} holds for all $k \in \N{}$.  Taking total expectations, it follows that, for all $k \in \N{}$,
  \bequationn
    \E[f(x_{k+1})] - \E[f(x_k)] \leq -\alpha \theta_1 \E[\|\nabla f(x_k)\|^2] + \theta_2.
  \eequationn
  Summing both sides for $k \in \{1,\dots,K\}$ yields
  \bequationn
    f_* - f(x_1) \leq \E[f(x_{K+1})] - f(x_1) \leq -\alpha \theta_1 \sum_{k=1}^K \E[\|\nabla f(x_k)\|^2] + K\theta_2.
  \eequationn
  Rearranging yields \eqref{eq.nonconvex_1}, then dividing by $K$ yields \eqref{eq.nonconvex_2}. 
\endproof

As in the case of \cite[Thm.~4.8]{BottCurtNoce18}, this result shows that while one cannot bound the expected optimality gap as when the P-L condition holds, one can bound the average norm of the gradients of the objective that are observed during the optimization process.

\subsection{No P-L Condition and Sublinearly Diminishing Stepsizes}

Finally, let us consider the behavior of \ref{alg.trish} when the P-L condition does not hold and diminishing stepsizes are employed.  For brevity, the following theorem considers both when parameters are chosen as in Theorem~\ref{th.diminishing_stepsize} and as in Theorem~\ref{th.diminishing_variance_diminishing_stepsize}, since in either case the final conclusion is the same.

\begin{theorem}\label{th.diminishing_variance_diminishing_stepsize_no_PL}
  Suppose Assumptions~\ref{ass.f} and \ref{ass.f_limits} hold and at least one of the following.
  \benumerate
    \item[(i)] With a pair $(h_1,h_2)$ satisfying the inequalities in Lemma~\ref{lem.fg_upperbound}, suppose that \ref{alg.trish} is run with $\{\gamma_{1,k}\}$, $\{\gamma_{2,k}\}$, and $\{\alpha_k\}$ chosen as in Theorem~\ref{th.diminishing_stepsize}.
    \item[(ii)] Suppose Assumption~\ref{ass.fg_upperbound2} holds and, with a pair $(h_3,h_4)$ satisfying the inequalities in Lemma~\ref{lem.fg_upperbound2}, suppose that \ref{alg.trish} is run with $\{\gamma_{1,k}\}$, $\{\gamma_{2,k}\}$, and $\{\alpha_k\}$ chosen as in Theorem~\ref{th.diminishing_variance_diminishing_stepsize}.
  \eenumerate
  Then, with $A_K := \sum_{k=1}^K \alpha_k$, it follows that
  \bsubequations\label{eq.nonconvex2}
    \begin{align}
      \lim_{K\to\infty} \E\left[ \sum_{k=1}^K \alpha_k \|\nabla f(x_k)\|^2 \right] &< \infty \label{eq.nonconvex2_1} \\
      \text{and}\ \ \E\left[ \frac{1}{A_K} \sum_{k=1}^K \alpha_k \|\nabla f(x_k)\|^2 \right] &\xrightarrow{K\to\infty} 0. \label{eq.nonconvex2_2}
    \end{align}
  \esubequations
\end{theorem}
\proof{Proof.}
  First observe that, under the conditions of the theorem, specifically the conditions placed on the stepsize sequence $\{\alpha_k\}$ in Theorem~\ref{th.diminishing_stepsize} or Theorem~\ref{th.diminishing_variance_diminishing_stepsize}, it follows that
\bequation\label{eq.robbins_monro}
    \sum_{k=1}^\infty \alpha_k = \infty\ \ \text{and}\ \ \sum_{k=1}^\infty \alpha_k^2 < \infty.
  \eequation
  Second, following the proofs of Theorem~\ref{th.diminishing_stepsize} or Theorem~\ref{th.diminishing_variance_diminishing_stepsize}, it follows that under conditions (i) or (ii) one finds by taking total expectations in \eqref{eq.for_later_1} or \eqref{eq.for_later_2} that
  \bequationn
    \baligned
      \E[f(x_{k+1})] - \E[f(x_k)] &\leq -\thalf \gamma_1 \alpha_k \E[\|\nabla f(x_k)\|^2] + \delta\alpha_k^2 \\
      \text{or}\ \ \E[f(x_{k+1})] - \E[f(x_k)] &\leq -\beta_1 \alpha_k \E[\|\nabla f(x_k)\|^2] + \beta_2\alpha_k^2.
    \ealigned
  \eequationn
  Without loss of generality, let us assume that condition (ii) holds and the latter inequality above is satisfied.  (The proof is the same if condition (i) holds and the former inequality above is satisfied.)  Summing both sides for $k \in \{1,\dots,K\}$ yields
  \bequationn
    f_* - f(x_1) \leq \E[f(x_{K+1})] - f(x_1) \leq -\beta_1 \sum_{k=1}^K \alpha_k \E[\|\nabla f(x_k)\|^2] + \beta_2 \sum_{k=1}^K \alpha_k^2,
  \eequationn
  which after rearrangement gives
  \bequationn
    \sum_{k=1}^K \alpha_k \E[\|\nabla f(x_k)\|^2] \leq \frac{f(x_1) - f_*}{\beta_1} + \frac{\beta_2}{\beta_1} \sum_{k=1}^K \alpha_k^2.
  \eequationn
  From \eqref{eq.robbins_monro}, it follows that the right-hand side converges to a finite limit as $K \to \infty$, giving~\eqref{eq.nonconvex2_1}.  Then, the limit \eqref{eq.nonconvex2_2} follows since \eqref{eq.robbins_monro} ensures that $\{A_K\} \to \infty$ as $K \to \infty$.
\endproof

A consequence of this theorem is the straightforward fact that
\bequationn
  \liminf_{k\to\infty} \E[\|\nabla f(x_k)\|^2] = 0.
\eequationn
That is, under the conditions of the theorem, the expected squared norms of the gradients at the iterates of the algorithm cannot stay bounded away from zero.

\section{Numerical Experiments}\label{sec.numerical}

In this section, we provide the results of numerical experiments to demonstrate the performance of \ref{alg.trish} compared to a stochastic gradient (SG) approach.  Through solving machine learning test problems involving objective functions of the form \eqref{prob.f_sum}---some convex and some nonconvex---we demonstrate that \ref{alg.trish} can outperform SG with comparable computational effort.  Before presenting our results, we first discuss how the parameters of the algorithm might be chosen.

\subsection{Algorithm Parameter Selection}\label{subsec.parameter}

Our analysis in \S\ref{sec.convergence} provides guidelines on how the stepsizes $\{\alpha_k\}$ and pairs $\{(\gamma_{1,k},\gamma_{2,k})\}$ should be chosen to guarantee convergence properties for \ref{alg.trish}.  That said, as for SG, the values required by the theory are often too conservative in practice, whereas one often finds better performance by a parameter tuning scheme.  Still, it is worthwhile to comment on how the theoretical analysis might inform parameter selection.  For our purposes, since our numerical experiments focus on results obtained with fixed parameters, we shall discuss how the analysis in~\S\ref{subsec.first} informs parameter selection.  Similar conclusions can be drawn based on our other theoretical results.

For simplicity, let us assume that the bound \eqref{eq.second_moment} in Assumption~\ref{ass.f_limits} holds with $M_2 = 1$.  In this case, the bound \eqref{eq.second_moment} is equivalent to the restriction that the variance of the stochastic gradient estimate is bounded by $M_1$, i.e., that $\E_k[\|g_k\|^2] - \|\nabla f(x_k)\|^2 \leq M_1$.  If one has an estimate $\widetilde{M}_1$ of $M_1$---which, for example, can be obtained by sampling gradients and computing a variance estimate---then, following Lemma~\ref{lem.fg_upperbound}, one can employ the value $\widetilde{h}_2 = \thalf\sqrt{\widetilde{M}_1} + 1$ for parameter selection.  In particular, Theorem~\ref{th.bounded_variance_fixed_stepsize} suggests to choose $(\gamma_1,\gamma_2)$ such that
\bequationn
  \frac{\gamma_1}{\gamma_2} < \frac{\widetilde{h}_2}{\widetilde{h}_2 - 1} = 1 + \frac{2}{\sqrt{\widetilde{M}_1}}.
\eequationn
Naturally, this still leads to flexibility in the precise values of $(\gamma_1,\gamma_2)$, but the trade-offs between different choices become similar to the traditional trade-offs one finds for the selection of $\alpha$ in an SG scheme: (i) one can choose values such that $\gamma_1-\gamma_2$ is large, which leads to fast convergence, but only to a relatively large neighborhood of the solution, or (ii) one can choose values such that $\gamma_1-\gamma_2$ is small, which leads to slow convergence, but to a relatively small neighborhood of a solution.  Overall, one might be discouraged by the idea that the choice of $(\gamma_1,\gamma_2)$ requires estimation of the upper bound~$M_1$.  However, this is not dissimilar to the fact that, theoretically, one needs an estimate of the Lipschitz constant $L$ of the gradient in order to choose the stepsize for SG, and clearly also for \ref{alg.trish}, such as through the bound \eqref{eq.fixed_alpha_choice}.  The good news is that estimating the variance of the stochastic gradient estimates is a reasonable request that could even be done during an initial phase that simply uses SG iterations.

Despite all of this commentary, in practice one should expect to achieve better performance by simply tuning parameters for a given problem, as is often done for SG methods.  For our experiments described in the following subsections, we chose $(\gamma_1,\gamma_2)$ by a simple tuning scheme that also selects the stepsize $\alpha$.  We took care to make sure that the tuning procedure for \ref{alg.trish} did not require more effort than the tuning used for the SG method that we have for comparison purposes.

\subsection{Logistic Regression}\label{sec.logistic_regression}

As a first test case, we considered the problem of binary classification through logistic regression using a few datasets available in the well-known LIBSVM repository; see~\cite{CC01a}.  In particular, for each dataset, with training feature vector $z_i \in \R{n}$ and training label $y_i \in \{-1,1\}$ for all $i \in \{1,\dots,N\}$, the objective of this problem has the form
\bequation\label{eq.logistic_regression}
  f(x) = \frac{1}{N} \sum_{i=1}^N \log(1 + e^{-y_i(x^Tz_i)}).
\eequation
Also available in each case is a testing dataset $\{(\zbar_i,\ybar_i)\}_{i=1}^{\Nbar}$.

We ran implementations of \ref{alg.trish} and SG and compare performance by comparing \emph{training loss} (i.e., the objective function~\eqref{eq.logistic_regression} evaluated with the training data) and \emph{testing accuracy} (i.e., for a given approximate solution, what fraction of the testing set is classified correctly) for iterates throughout the optimization process.  We ran each algorithm for one epoch (i.e., until $N$ training pairs have been accessed) with a fixed stepsize $\alpha$ and, for \ref{alg.trish}, a fixed parameter pair $(\gamma_1,\gamma_2)$.

For both algorithms and all datasets, the stochastic gradient estimates were computed using a mini-batch size of 64.  For choosing a fair set of parameters for the comparison for each dataset, we first ran SG with a stepsize of $0.1$ and computed $G$ as the average norm of stochastic gradient estimates throughout the run.  Then, for~\ref{alg.trish}, we considered the stepsizes $\alpha \in \{10^{-1},10^{-1/2},10^0,10^{1/2},10^1\}$ and parameters $\gamma_1 \in \{\tfrac4G, \tfrac8G, \tfrac{16}G, \tfrac{32}G\}$ and $\gamma_2 \in \{\tfrac1{2G},\tfrac1G,\tfrac2G\}$.  (The value~$G$ gauges the magnitude of the stochastic gradient estimates, which depends on problem scaling.  As seen in our results, these choices of $(\gamma_1,\gamma_2)$ ensure that step normalization---i.e., case~2 of \ref{alg.trish}---occurs.  In practice, one could compute~$G$ during an initial SG phase before starting \ref{alg.trish}, but to cleanly distinguish between \ref{alg.trish} and SG, we computed this value using an independent run of SG.)  This resulted in 60 parameter settings with \ref{alg.trish} employing stepsizes in the range from $\tfrac1{2G} \times 10^{-1}$ (i.e., the minimum $\gamma_2$ times the minimum $\alpha$) to $\tfrac{32}G \times 10^1$ (i.e., the maximum $\gamma_1$ times the maximum $\alpha$).  Hence, for SG, we considered 60 values for $\alpha$ in the range $[\tfrac1{2G} \times 10^{-1},\tfrac{32}G \times 10^1]$ so that neither algorithm had an advantage in terms of the range of the stepsizes.  Specifically, we considered the~60 values such that $\log_{10}(\alpha)$ was evenly distributed in $[\log_{10}(\tfrac1{2G} \times 10^{-1}),\log_{10}(\tfrac{32}G \times 10^1)]$.

For each dataset, we ran the algorithms with these different parameters settings and selected for each the setting that led to the best average testing accuracy in the last ten iterations of the run.

\subsubsection{a1a.}

The first dataset that we considered was \texttt{a1a} in which the feature vectors have length $n = 123$, the number of points in the training set is $N = 1605$, and the number of points in the testing set is $\Nbar = 30956$.  For tuning, the value $G \approx 0.1746$ was determined, yielding a stepsize range of approximately $[0.2863,1832]$.  After tuning, the selected parameter setting for \ref{alg.trish} was $(\alpha,\gamma_1,\gamma_2) \approx (0.1,22.90,2.863)$ and the selected parameter setting for SG was $\alpha \approx 0.4471$.

The algorithms, \ref{alg.trish} and SG, were each run 10 times from the same starting point (the origin).  The training losses and testing accuracies, averaged over these 10 runs, are plotted in Figure~\ref{fig.a1a} after 0.1 epoch through the end of the first epoch.  (The values during the first 0.1 epoch are not plotted here, nor for the other datasets, so that it is easier to distinguish the differences at the end of the first epoch.)  It is worthwhile to note that, during the runs for \ref{alg.trish}, case 1 did not occur, case 2 occurred in approximately 99\% of the iterations, and case 3 occurred in approximately 1\% of the iterations; i.e., step normalization occurred in a large majority of the iterations.  The figure shows that \ref{alg.trish} yields better training losses throughout the optimization process.  However, for this dataset, the performance in terms of testing accuracy is roughly the same for both algorithms.

\bfigure[ht]
  \centering
  \begin{subfigure}{.49\textwidth}
    \centering
    \includegraphics[width=1\linewidth]{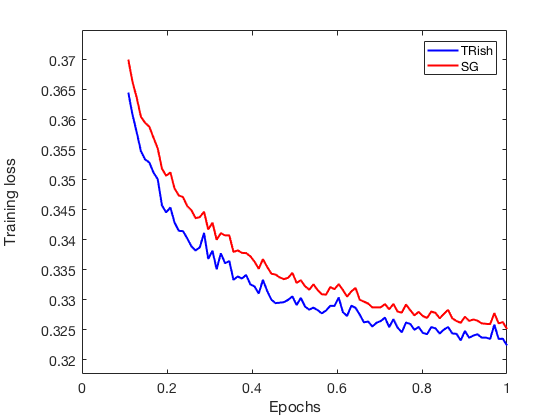}
  \end{subfigure}
  \begin{subfigure}{.49\textwidth}
    \centering
    \includegraphics[width=1\linewidth]{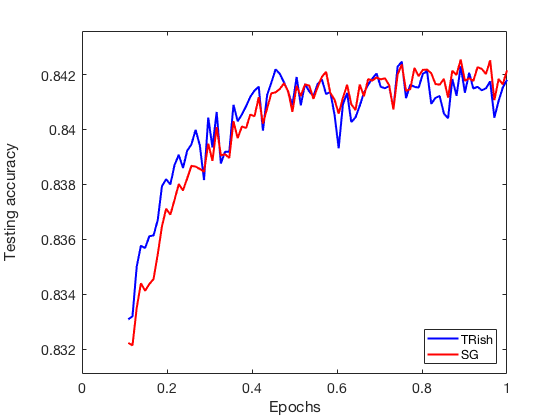}
  \end{subfigure}
  \caption{Average training loss and testing accuracy during the first epoch when \ref{alg.trish} and SG are employed to minimize the logistic regression function \eqref{eq.logistic_regression} using the \texttt{a1a} dataset.}
  \label{fig.a1a}
\efigure

\subsubsection{w1a.}

The second dataset that we considered was \texttt{w1a} in which the feature vectors have length $n = 300$, the number of points in the training set is $N = 2477$, and the number of points in the testing set is $\Nbar = 47272$.  For tuning, the value $G \approx 0.0887$ was determined, yielding a stepsize range of approximately $[0.5638,3608]$.  After tuning, the selected parameter setting for \ref{alg.trish} was $(\alpha,\gamma_1,\gamma_2) \approx (0.1,360.8,5.638)$ and the selected parameter setting for SG was $\alpha \approx 0.6541$.

The training losses and testing accuracies, averaged over 10 runs when both algorithms were initialized at the same starting point (the origin), are plotted in Figure~\ref{fig.w1a}.  During the runs for \ref{alg.trish}, case 2 occurred in approximately 99\% of the iterations while case 1 and case 3 combined occurred in fewer than 1\% of the iterations.  For this dataset, \ref{alg.trish} outperformed SG both in terms of training losses and testing accuracies throughout the first epoch.

\bfigure[ht]
  \centering
  \begin{subfigure}{.49\textwidth}
    \centering
    \includegraphics[width=1\linewidth]{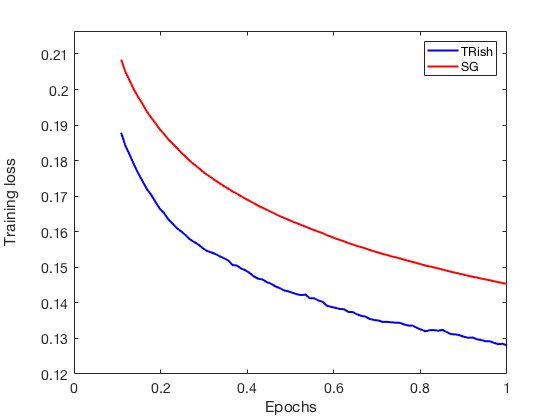}
  \end{subfigure}
  \begin{subfigure}{.49\textwidth}
    \centering
    \includegraphics[width=1\linewidth]{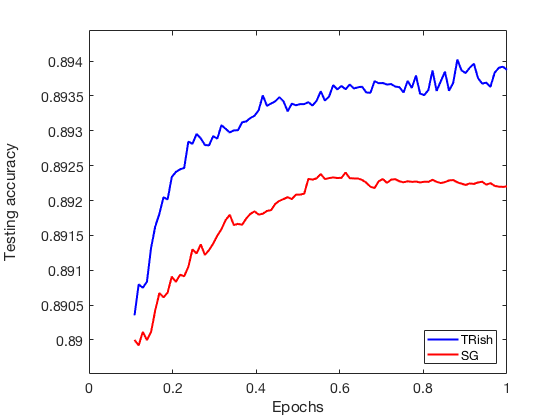}
  \end{subfigure}
  \caption{Average training loss and testing accuracy during the first epoch when \ref{alg.trish} and SG are employed to minimize the logistic regression function \eqref{eq.logistic_regression} using the \texttt{w1a} dataset.}
  \label{fig.w1a}
\efigure

\subsubsection{rcv1.}

The third dataset that we considered was \texttt{rcv1} in which the feature vectors have length $n = 47236$, the number of points in the training set is $N = 20242$, and the number of points in the testing set is $\Nbar = 677399$.  For tuning, the value $G \approx 0.0497$ was determined, yielding a stepsize range of approximately $[1.007,6444]$.  After tuning, the selected parameter setting for \ref{alg.trish} was $(\alpha,\gamma_1,\gamma_2) \approx (0.3162,644.4,10.07)$ and the selected parameter setting for SG was $\alpha \approx 10.84$.

The training losses and testing accuracies, averaged over 10 runs when both algorithms were initialized at the same starting point (the origin), are plotted in Figure~\ref{fig.rcv1}.  During the runs for \ref{alg.trish}, case 1 occurred in approximately 27\% of the iterations, case 2 occurred in approximately~73\% of the iterations, and case 3 did not occur.  For this dataset, \ref{alg.trish} outperformed SG both in terms of training losses and testing accuracies throughout the first epoch.  That said, the testing accuracies appear to near at the end of the first epoch, leading one to wonder about the performance of the methods if the parameters are re-tuned and the algorithms are run for more epochs.

\bfigure[ht]
  \centering
  \begin{subfigure}{.49\textwidth}
    \centering
    \includegraphics[width=1\linewidth]{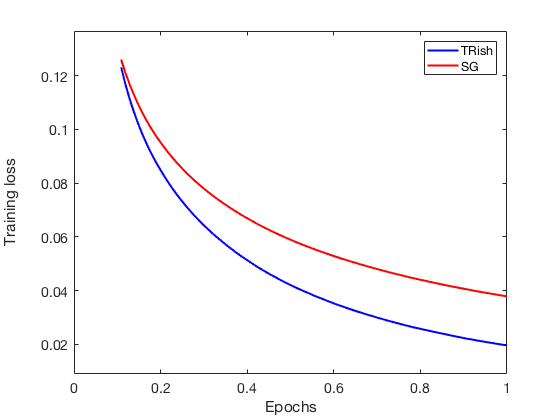}
  \end{subfigure}
  \begin{subfigure}{.49\textwidth}
    \centering
    \includegraphics[width=1\linewidth]{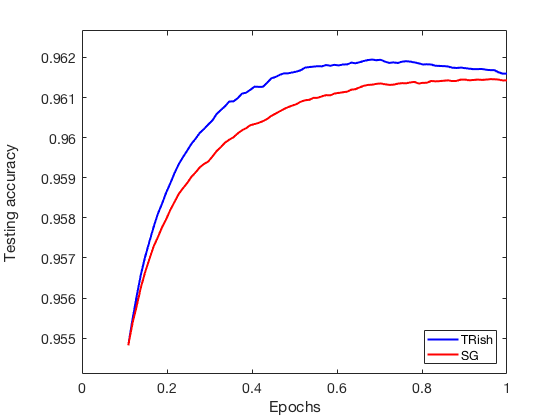}
  \end{subfigure}
  \caption{Average training loss and testing accuracy during the first epoch when \ref{alg.trish} and SG are employed to minimize the logistic regression function \eqref{eq.logistic_regression} using the \texttt{rcv1} dataset.}
  \label{fig.rcv1}
\efigure

To address this question, Figure~\ref{fig.rcv1_2} plots the training losses and testing accuracies---averaged over 10 runs---for \ref{alg.trish} and SG during two epochs.  (For this horizon, tuning led to the parameter setting for \ref{alg.trish} as $(\alpha,\gamma_1,\gamma_2) = (0.1,376.2,47.02)$ and the parameter setting for SG as $\alpha \approx 5.192$.  For \ref{alg.trish}, case~2 occurred in approximately 94\% of the iterations, case~3 occurred in approximately~5\% of the iterations, and case~1 occurred in fewer than 1\% of the iterations.)  These plots show a trade-off where, for a longer horizon, the better parameters for \ref{alg.trish} do not necessarily offer better results initially, but do offer better results eventually.

\bfigure[ht]
  \centering
  \begin{subfigure}{.49\textwidth}
    \centering
    \includegraphics[width=1\linewidth]{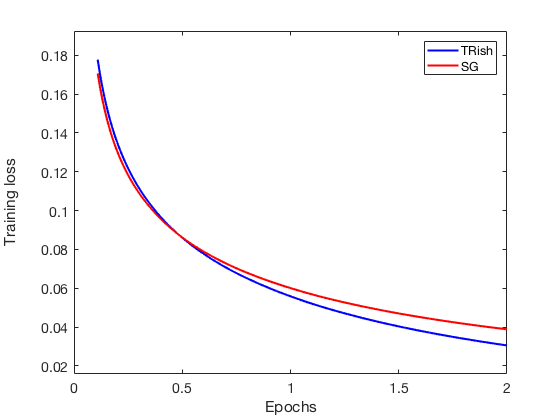}
  \end{subfigure}
  \begin{subfigure}{.49\textwidth}
    \centering
    \includegraphics[width=1\linewidth]{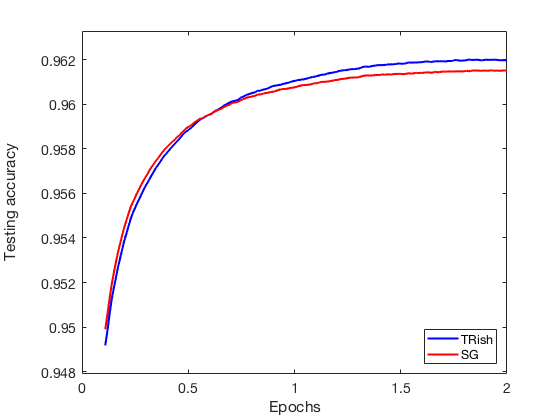}
  \end{subfigure}
  \caption{Average training loss and testing accuracy during the first two epochs when \ref{alg.trish} and SG are employed to minimize the logistic regression function \eqref{eq.logistic_regression} using the \texttt{rcv1} dataset.}
  \label{fig.rcv1_2}
\efigure

In all of the experiments presented in this section, \ref{alg.trish} generally outperforms SG.  However, the gains are somewhat limited due to the fact that, by convexity of the problems, both algorithms are tending to neighborhoods around the same optimal solution.  The results presented in the next subsection, in which we consider nonconvex optimization problems arising from neural network training, show more substantial benefits from using \ref{alg.trish} as compared to SG.

\subsection{Neural Network Training}

As a second test case, we considered the problem to train convolutional neural networks (CNNs) for image classification.  We considered two well-known datasets.  The first, the \texttt{mnist} dataset \cite{LeCuBottBengHaff98}, is a collection of images of hand-written digits.  The goal for training the network for this dataset is to classify which of the digits (0 through 9) is written in each image.  It includes $N = 60000$ training samples and $\Nbar = 10000$ testing samples.  The second, the \texttt{cifar-10} dataset \cite{CIFAR10}, is a collection of color images in ten categories (e.g., airplanes, dogs, and ships).  The goal for training the network for this dataset is to classify the image with the correct category.  It includes $N = 50000$ training samples and $\Nbar = 10000$ testing samples.

Implemented using \texttt{tensorflow}, the neural networks that we considered for both datasets are composed of two convolutional layers (involving 32 and 64 filters, respectively, and each followed by an average pooling layer) followed by two fully connected layers.  ReLU activation is used at each hidden layer and the objective is defined using the logistic (cross entropy) loss function.  The networks vary slightly, e.g., due to the fact that a pixel for each \texttt{mnist} image corresponds to a single feature while a pixel for each \texttt{cifar-10} image corresponds to three features (for each RGB value since they are color images).  As seen in our experimental results, training the network led to a very good classifier for \texttt{mnist}, yielding over 95\% testing accuracy.  The performance is less impressive for \texttt{cifar-10} (yielding around 60\% accuracy); achieving higher accuracy would require a more sophisticated network and more computational resources than were available.  That said, both datasets provide interesting settings for comparing the performance of \ref{alg.trish} and SG.

As for the results in \S\ref{sec.logistic_regression}, we compare performance between \ref{alg.trish} and SG by comparing training loss and testing accuracy.  We tuned parameters using the same setup as in \S\ref{sec.logistic_regression}, except with slightly different parameter choices.  In particular, the mini-batch size used when computing stochastic gradients was 128 and, when computing $G$, we ran SG with a stepsize of 0.01.  For \ref{alg.trish}, we considered stepsizes $\alpha \in \{10^{-3},10^{-2},10^{-1},10^0\}$ and parameters $\gamma_1 \in \{\tfrac4G,\tfrac8G,\tfrac{16}G\}$ and $\gamma_2 \in \{\tfrac1{8G},\tfrac1{4G},\tfrac1{2G}\}$.  This means that SG was tuned with 36 choices of $\alpha$ in the range $[\tfrac1{8G}\times10^{-3},\tfrac{16}G\times10^0]$.

\subsubsection{\texttt{mnist}.}

For \texttt{mnist}, we ran the algorithms for two epochs.  For parameter tuning, the value $G \approx 2.8277$ was determined, yielding a stepsize range of approximately $[2.683\times10^{-5},3.435]$.  After tuning, the selected parameter setting for \ref{alg.trish} was $(\alpha,\gamma_1,\gamma_2) \approx (1,1.717,0.0268)$ and the selected parameter setting for SG was $\alpha \approx 0.0609$.

The training losses and testing accuracies for each of the 10 runs that we performed with the tuned parameters are plotted in Figure~\ref{fig.mnist}, ignoring the first 0.2 epochs so that the later values are more easily distinguished.  (For each run, the network parameters were initialized to the same randomly generated values; the values were generated from a truncated normal distribution with mean 0 and standard deviation $0.1$.  We did not average the loss and accuracy values over the 10 runs since the optimization problem is nonconvex, meaning that for each run an algorithm might tend toward a different region of the search space.)  During the runs for \ref{alg.trish}, case 1 occurred in approximately 62\% of the iterations, case 2 occurred in approximately~37\% of the iterations, and case 3 almost did not occur.  Overall, \ref{alg.trish} consistently outperformed SG in terms of both training loss and testing accuracy throughout the optimization process.

\bfigure[ht]
  \centering
  \begin{subfigure}{.49\textwidth}
    \centering
    \includegraphics[width=1\linewidth]{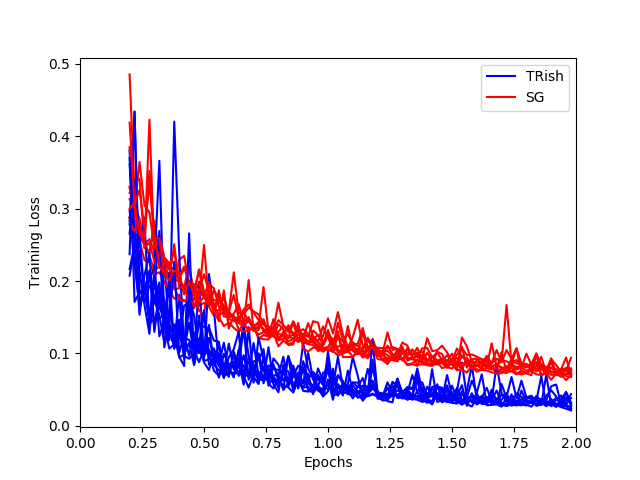}
  \end{subfigure}
  \begin{subfigure}{.49\textwidth}
    \centering
    \includegraphics[width=1\linewidth]{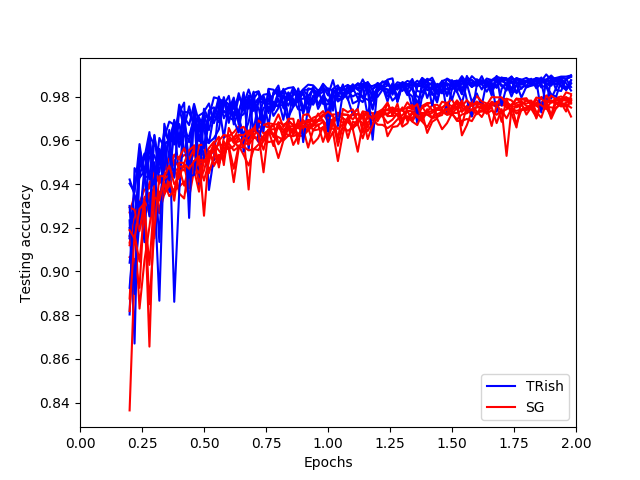}
  \end{subfigure}
  \caption{Average training loss and testing accuracy during the first two epochs when \ref{alg.trish} and SG are employed to train a convolutional neural network using the \texttt{mnist} dataset.}
  \label{fig.mnist}
\efigure

\subsubsection{\texttt{cifar-10}.}

For \texttt{cifar-10}, we ran the algorithms for five epochs (since further improvement was clearly being made even after the first few epochs).  The value $G \approx 964.39$ was determined, yielding a stepsize range of approximately $[8.990\times10^{-6},1.151]$.  After tuning, the parameter setting for \ref{alg.trish} was $(\alpha,\gamma_1,\gamma_2) \approx (1,1.051,0.0089)$ and the parameter setting for SG was $\alpha \approx 0.0104$.

The training losses and testing accuracies for each of the 10 runs that we performed with the tuned parameters are plotted in Figure~\ref{fig.cifar10}, again ignoring the first 10\% of the runs (i.e., in this case, the first 0.5 epochs) so that the later values are more easily distinguished.  (For each run, the network parameters were initialized to the same randomly generated values; the values were generated from a truncated normal distribution with mean 0 and standard deviation $0.01$.)  During the runs for \ref{alg.trish}, case 1 occurred in approximately 1\% of the iterations, case 2 occurred in approximately~99\% of the iterations, and case 3 did not occur.  In these experiments, \ref{alg.trish} typically outperformed SG in terms of both training loss and testing accuracy throughout each run.

\bfigure[ht]
  \centering
  \begin{subfigure}{.49\textwidth}
    \centering
    \includegraphics[width=1\linewidth]{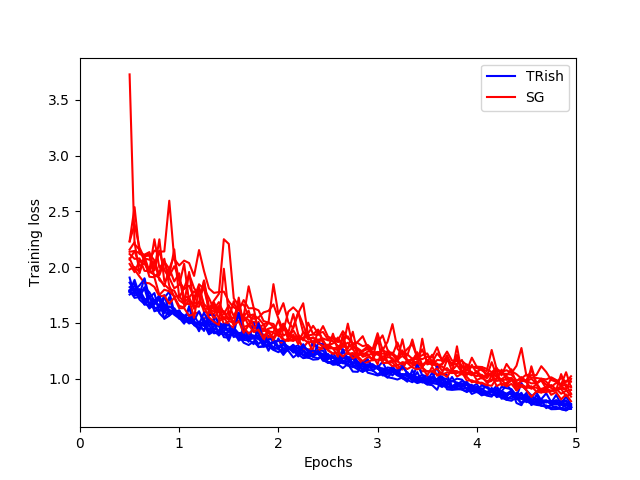}
  \end{subfigure}
  \begin{subfigure}{.49\textwidth}
    \centering
    \includegraphics[width=1\linewidth]{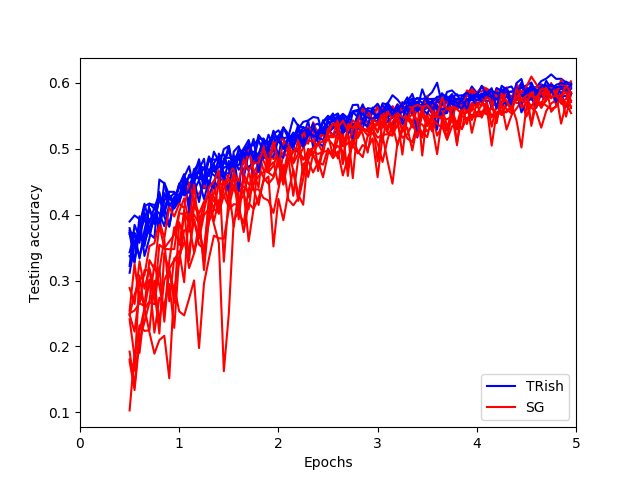}
  \end{subfigure}
  \caption{Average training loss and testing accuracy during the first five epochs when \ref{alg.trish} and SG are employed to optimize the convolutional neural network using the \texttt{cifar10} dataset.}
  \label{fig.cifar10}
\efigure

\section{Conclusion}\label{sec.conclusion}

An algorithm inspired by a trust region methodology has been proposed, analyzed, and tested for solving stochastic and finite-sum minimization problems.  Our proved theoretical guarantees show that our method, deemed \ref{alg.trish}, has convergence properties that are similar to a traditional SG method.  Our numerical results, on the other hand, show that \ref{alg.trish} can outperform a traditional SG approach.  We attribute this better behavior to the algorithm's use of normalized steps, which one can argue lessens its dependence on problem-specific quantities.

Naturally, a more substantial numerical study---that goes well beyond the scope of this paper---would be necessary to fully explore the trade-offs between \ref{alg.trish} and SG in practice.  For example, a more substantial numerical study would take into account different procedures that might be used to decrease the stepsize after some number of iterations, as is typically done in practice.  Indeed, for the convex problems that we considered, this was our motivation for presenting results for only one epoch, since, in practice, one often adjusts the stepsize after each epoch.  For \ref{alg.trish}, this adjustment may involve updates to the pair $(\gamma_1,\gamma_2)$ as well, which one might adjust so that $\gamma_1 - \gamma_2 = \Ocal(\alpha)$, as our theory suggests.

Finally, while not considered in this paper, we believe it would be interesting to explore the incorporation within \ref{alg.trish} of various enhancements, such as the use of second-derivative (i.e., Hessian) approximations, acceleration, and/or momentum.  These might further improve the practical performance of the framework set forth in this paper.

\ifthenelse{\coralreport = 1}{\section*{Acknowledgment} The authors would like to thank Chaoxu Zhou of Columbia University for his valuable assistance in correcting errors in an earlier version of this manuscript.  They would also like to thank the anonymous referees and the Associate Editor for providing valuable comments on an earlier draft of the paper.  In particular, one anonymous referee suggested a proof for Lemma~\ref{lem.fg_upperbound}, which greatly improved the analysis.  All authors were supported by the U.S.~National Science Foundation's Division of Computing and Communication Foundations and Division of Mathematical Sciences under award numbers CCF--1618717 and DMS--1319356.  The first author was also supported by the U.S.~Department of Energy, Office of Science, Applied Mathematics, Early Career Research Program under Award Number DE--SC0010615.}{\ACKNOWLEDGMENT{The authors would like to thank Chaoxu Zhou of Columbia University for his valuable assistance in correcting errors in an earlier version of this manuscript.  They would also like to thank the anonymous referees and the Associate Editor for providing valuable comments on an earlier draft of the paper.  In particular, one anonymous referee suggested a proof for Lemma~\ref{lem.fg_upperbound}, which greatly improved the analysis.  All authors were supported by the U.S.~National Science Foundation's Division of Computing and Communication Foundations and Division of Mathematical Sciences under award numbers CCF--1618717 and DMS--1319356.  The first author was also supported by the U.S.~Department of Energy, Office of Science, Applied Mathematics, Early Career Research Program under Award Number DE--SC0010615.}}

\ifthenelse{\coralreport = 1}{
\bibliographystyle{plain}
}{
\bibliographystyle{informs2014}
}
\bibliography{references}

\end{document}